\documentclass[12pt]{article}
\usepackage{amssymb,latexsym,amsfonts}
\def\R{{\mathbb R}}
\def\Zp{{\mathbb Z}_+}
\def\C{\mathbb C}
\def\Z{\mathbb Z}

\def\Hi{\ell^2(\Zp)}

\def\da{a^{\dag}}

\def\db{b^{\dag}}
\def\A{${\cal A}_q$ }
\def\al{a(\lambda)}
\def\dal{a^{\dag}(\lambda)}
\def\H0{${\cal H}_0$ }
\def\Hg{${\cal H}_{\gamma}$ }
\def\ba{\begin{array}}
\def\ea{\end{array}}
\def\l{\label}
\def\ra{\rightarrow}
\def\qD{{}^q\!D_z}
\def\ba{\begin{array}}
\def\ea{\end{array}}
\def\l{\label}
\def\ds{\displaystyle}
\newcommand{\beq}{\begin{equation}}
\newcommand{\eeq}{\end{equation}}
\newcommand{\be}{\begin{equation}}
\newcommand{\ee}{\end{equation}}
\begin{document}
\begin{titlepage}
\thispagestyle{empty}
\begin{flushright}
\vspace*{.5cm}
{\large SPBU-IP-97-31
\\ math.QA/9803xxx}
\end{flushright}
\vspace*{.5cm}
\begin{center}
{\Large \bf On position operator spectral measure}\\[.5cm]
{\Large \bf for deformed oscillator in the case}\\[.5cm]
{\Large \bf of indetermine Hamburger moment problem
\footnote{This research was supported by RFFI
grants No 97-01-01152 (VVB, EVD) and No 96-01-00851 (PPK).}}
\end{center}
\vskip .3cm

\centerline{\Large \bf V.V.Borzov${}^{\left.1\right)}$,\
E.V.Damaskinsky${}^{\left.2\right)}$\
and P.P.Kulish${}^{\left.3\right)}$}
\begin{center}
{${}^{\left.1\right)}$\ Sanct-Petersburg University of
Telecommunications, St.Petersburg, Russia}
\\
{${}^{\left.2\right)}$\
University of Defense Technical
Engineering, St.Petersburg, Russia}\\
{${}^{\left.3\right)}$\
Department of the Steklov Mathematical Institute,
St.Petersburg,
Russia}
\end{center}

\vskip 0.5cm

\begin{abstract}
The spectral measure of the position (momentum) operator $X$
for $q$-deformed oscillator is calculated
in the case of the indetermine Hamburger moment problem.
The exposition is given for concrete choice of
generators for $q$-oscillator algebra,
although developed technique
apply for every other cases with indetermine moment problem.
The Stieltjes transformation $m(z)$ of spectral measure is
expressed in terms of the entries of Jacobi matrix $X$ only.
 The direct connection between values of parameters
labeling the spectral measures and related selfadjoint
extensions of $X$ is established.
\end{abstract}
\end{titlepage}
\tableofcontents

\section{Introduction}
\renewcommand{\theequation}{\thesection.\arabic{equation}}
\setcounter{equation}{0}

The connection of the quantum harmonic oscillator with the
Hermite polynomials $H_n(x)$ is well-known from the early days
of quantum theory. The coordinate operator
$X=\ds{\frac{1}{\sqrt{2}}}\left( b + \db \right)$,
being the Jacobi matrix in the
number representation
$$
N|n\rangle = \db b |n\rangle = n |n\rangle\, , \qquad
n=0,1,2,\ldots\, ,
$$
has the eigenvectors, components of  which in coordinate
representation are proportional to $H_n(x)$.

The discovery of the quantum groups gave rise to connection of
quantum algebras to variety of $q$-special functions. In
particular, it was pointed out that deformed oscillator algebra
\A \cite{L1}-\cite{L6}
is related to $q$-Hermite polynomials
(see e.g. \cite{L7,L8,Z199}). The possibility of different
choice of $q$-oscillator generators $\al$, $\dal$
results in \cite{Z189}-\cite{BDY}
different behaviour of the corresponding Jacobi matrix entries
$$
\left( X(q,\lambda)\right)_{i,j}=b_i{\delta}_{i,j-1}
+ b_{i-1}{\delta}_{i,j+1}
$$
and different $q$-Hermite polynomials $H_n(x;q,\lambda)$.

The spectral measure of a selfadjoint extension
of the symmetric operator $X(q,\lambda)$.
and of the orthogonal $q$-polynomials is
related to the classical
Hamburger moment problem\cite{L10}-\cite{Chihara2}
This moment problem can be determined (the measure is unique)
or indetermine (a family of measures) according to the
behaviour of the Jacobi matrix entries \cite{L10,L11}. It
was pointed out in \cite{Nagel} that even for the quantum
harmonic oscillator the Hamburger moment problem for the
Jacobi matrices $J^{(k)}=b^k + ({\db})^k$ (appearing in the
description of higher power squeezed states) is indetermine
for $k>2$.  As far as we know, for all the choices of
$q$-oscillator generators \cite{L13,L14}, when the
corresponding $q$-Hermite polynomials and the spectral
measures where known from the $q$-analysis
\cite{Ex,Gasper,Andrews} the moment problems are determined.

The calculation of the spectral measure of the position
(momentum) operator $X(q,\lambda)$ in the case of the
indetermine Hamburger moment problem is the aim of this paper.
Although the developed construction applies to each choice
of generators of the $q$-oscillator
algebra \A
the particular expressions and proofs are given for the
$q$-oscillator \cite{L4,L5} with nonzero entries of
the Jacobi matrix equal to
$b_n=\sqrt{[n+1]}$,
where "symmetric" basic
number $[a]$ being defined as
$[a]=\ds{\frac{q^a-q^{-a}}{q-q^{-1}}}.$
Let us note that $[ a ]$ is exponentially
growing for $q>1,$ as well as  for $0<q<1$ in contrast with the
$q$-number $[a;q]=\ds{\frac{1-q^{a}}{1-q}}$
usual for $q$-analysis \cite{Andrews} and basic
hypergeometric functions \cite{Ex,Gasper}.

Following to the general theory of the moment problem and
Jacobi matrices \cite{L10,L11,L12} we express the Stieltjes
transformation $m(z)$ of the spectral measure in terms of the
entries $b_n$ of Jacobi matrix only. We also establish
the direct connection between the values of parameters
labeling the spectral measures and associated selfadjoint
extensions of position operator $X$, which was missing in the
general considerations.

The paper is organized as follows. The relations
between harmonic oscillator and Hermite
polynomials, main formulas of the moment
problem and definitions of the $q$-deformed
oscillator algebra \A and $q$-Hermite polynomials
\cite{Z199} are briefly reviewed in Sec.2. In this section
the Hamburger power moment problem for deformed
oscillator of considered type is also formulated.
This moment problem is indetermine one. In this case
extremal spectral measure is concentrated
\cite{L10,L11,L12} on the set of zeros of entire
function expressed in terms of the related
orthogonal polynomials. Let us note that in the cases
considered, for example, in the papers
\cite{Chihara,b1,b4} where such measures, also
for indetermine Hamburger moment problem,
are constructed the founded measures are expressed in
terms of the standard for q-analysis symbol
$ [\alpha ;q]$
and represent the q-analogues of some classical special
functions.
The entries of the Jacobi matrix for the $q$-oscillator
generators under concideration are expressed in terms
of the symmetric $q$-number $[n]$ or its one parameter
generalization
$[n;q,\lambda] = q^{\lambda n}[n]$.
The spectral measures we are interested in, are
expressed  as infinite series with this $q$-numbers.
These series were not study yet in the $q$-analysis and they
could be of interest by themselves.
Our main considerations and computations are given in Sec.3.
In this section we first compute the value  $m(i)$
of the Stieltjes transformation $m(z)$ of the spectral measure
at the point $i$, and after that we found the value of
$m(z)$  at arbitrary complex value $z$.
Then in the Sec.4 we give the construction of the support of
spectral measure $\ds{\sigma _{\varphi _0}}$.
In Sec.5 concrete
examples of the spectral measures $\sigma _0$ and
$\sigma _\pi $ are given.

\section{Background material on q-oscillator,
classical \\ moment problem
and $q$-Hermite polynomials}
\renewcommand{\theequation}{\thesection.\arabic{equation}}
\setcounter{equation}{0}

\subsection{Spectral theory of Jacobi matrices,
classical moment problem and q-ortho\-go\-nal polynomials}

We  recall some of the results on the spectral
theory of Jacobi matrices and their relation with orthogonal
polynomials\footnote{For more information we refer to the
books by Berezanski\u\i\ \cite{L11} [Ch.~VII, \S 1],
Ahiezer\ \cite{L10} and Shohat and Tamarkin\ \cite{L12}}.

Let operator $X$ acting on the standard orthonormal basis
$\{ e_n\mid n\in\Zp\}$ of $\Hi$ by
\begin{equation}\label{ad1}
Xe_n =b_n\, e_{n+1} + a_n\, e_n + b_{n-1}\, e_{n-1},
\qquad a_n\geq 0,\ b_n\in{\R}.
\end{equation}
Then $X$ can be represented
(in  $\Hi$) by the
infinite dimensional three-diagonal matrix
\begin{equation}\label{ad1a}
X=\left(
\begin{array}{cccccc}
a_0 & b_0 &  0  &  0  & \cdots  & \cdots  \\
b_0 & a_1 & b_1 &  0  & \cdots  & \cdots  \\
 0  & b_1 & a_2 & b_2 & \cdots  & \cdots  \\
 0  & 0   & b_2 & a_3 & \ddots  & \cdots  \\
\vdots  & \vdots  & \vdots  & \ddots  & \ddots  & \ddots  \\
\vdots  & \vdots  & \vdots  & \vdots  & \ddots  & \ddots
\end{array}
\right) \, ,
\end{equation}
which is known as a Jacobi matrix.

With such Jacobi matrix one  can associate
polynomials $P_n(x)$ of degree $n$ in
$x$ by the recurrence relation
\begin{equation}\label{ad2}
b_n\, P_{n+1}(x)+a_n\, P_n(x)+b_{n-1}P_{n-1}(x) = xP_n(x)\, ,
\end{equation}
with   "initial conditions"
\begin{equation}\label{ad2a}
 P_0(x)=1\, ,\qquad P_{-1}(x)=0\, .
\end{equation}

In the following we assume that $a_n,\ b_n \in {\R}$ and
$a_n=0\ , b_n>0$.  Under this condition all polynomials $P_n(x)$
have real coefficients and $P_n(-x)= (-1)^n P_n(x).$

Note that the recurrence relation (\ref{ad2}) has
two linearly independent solutions.
The polynomials $P_n(x)$ solving the recurrence relation
(\ref{ad2}) and fulfilling the initial conditions (\ref{ad2a})
are called the polynomials of the  first kind.
One has ${\rm deg}\,P_n(x)=n$.

The independent set of solutions of the
relation (\ref{ad2}) consists of polynomials $Q_n(x)$
fulfilling the initial conditions
\begin{equation} \label{ad2b}
Q_0(x) =0,\qquad\qquad Q_1(x)=\frac{1}{b_0}\, ,
\end{equation}
and we have ${\rm deg}\,Q_n(x)=n-1$.
Such polynomials are called the polynomials of the
second kind for the Jacobi matrix $X$ (\ref{ad1a}).

These polynomials of the first and second
kind are related through
\begin{equation} \label{ad3c}
P_{n-1}(x)Q_{n}(x) - P_{n}(x)Q_{n-1}(x)=\frac{1}{b_{n-1}},
\qquad (n=1,2,3,...)\, .
\end{equation}

With the Jacobi matrix $X$ (\ref{ad1a}) one can
associate the so-called Hamburger (power) moment
problem \cite{L10}. For given number sequence
$s_n$ find the measure $\sigma $ on the line such that
\be\label{ad7}
s_n =\int_{-\infty }^\infty t^n d\sigma(t)\, ,
\quad n=0,1,2,\ldots\ ; \qquad
s_n =\alpha _0\, ,\qquad t^n=
\sum_{i=0}^n\alpha _iP_i(t;q)\, .
\ee
If such measure is defined uniquely then the
moment problem is determined one.
One says that the moment problem
is indetermine if one has infinite family of measures which
solves the relation (\ref{ad7}) for given sequence of moments
$s_n$.  There is one to one correspondence between $X$ and
$s_n$~\cite{L10,L11}~:
\be\l{k1}
s_n = \left( e_0, X^n e_0 \right)\, .
\ee

The matrix (\ref{ad1a}) defines the operator $X$ which
 is symmetric, and its deficiency indices
are $(0,0)\,$ (which mean that its closure $\overline{X}$ is a
selfadjoint operator)\, or $(1,1)\,$
(in this case $\overline{X}$
has many selfadjoint
extensions).

Note that the deficiency indices are $(0,0)$,
and thus $\overline{X}$
is selfadjoint in $\Hi$, if and only if
the related orthonormal polynomials $P_n(x)$,
correspond to a determined moment problem.
By Favard's theorem, for a selfadjoint operator $\overline{X}$
there exists a unique positive measure $m$ on the real line such
that the polynomials $P_n(x)$ are orthonormal
\begin{equation}\label{ad3}
\int_{\R} P_n(x)P_m(x) \, dm(x) = \delta_{n,m}.
\end{equation}
The measure is obtained by
$m(B)=\langle E(B)e_0,e_0\rangle$,
where $B\subset\R$ is a Borel set, and  $E$ denotes the
spectral decomposition of the selfadjoint operator $\overline{X}$.
The point spectrum corresponds to the
discrete mass points in $dm$. The spectral decomposition $E$ of
$\overline{X}$ is related to the orthogonality measure by
$$\langle E(B)e_n,e_m\rangle =
\int_B P_n(x)P_m(x) \, dm(x)\, ,$$ where
$B\subset\R$ is a Borel set.

If, in general case $a_n\neq 0$,
$|a_n|\leq C$, $n=0,1.2,\ldots$, and there exist $N,$
such that  the inequality $b_{n-1}b_{n+1}\leq b_n^{\, 2}$
is hold for all $n>N$ then in the case when
\begin{equation}
\label{ad3a}
\sum_{n=0}^\infty \frac 1{b_n}<\infty ,
\end{equation}
operator  $\overline{X}$
is not selfadjoint operator, and has infinitely many
selfadjoint extensions.

The operator  $\overline{X}$ has
deficiency indices $(0,0)$
$iff$ the series $\sum_{n=0}^\infty |P_n(z)|^2$
is divergent for
all $z\in \C$, ${\rm Im\,} z\neq 0$.
If the operator
$\overline{X}$ has  deficiency indices $(1,1)$ then this series
converges for all $z\in \C$, ${\rm Im\,} z\neq 0$.
 In the case when $\overline{X}$ has  deficiency indices $(1,1)$
the deficiency subspaces $N_{\overline{z}}$ are all
one-dimensional ones and spanned by the vectors
$\sum_{n=0}^\infty P_n(z)|n>$

The spectral properties of selfadjoint extensions
of operator $X$ are intimately linked with the
properties of measures which solves the related
moment problem. Let a parameter  $\varphi _0$
label \cite{a12} selfadjoint extensions
$X_{\varphi _0}$ of the
operator $X$ in the case of indetermine moment
problem. Then the spectral measure
$\sigma _{\varphi _0}$ gives the ''extremal''
solution of the Hamburger moment problem
(\ref{ad7}) related to the selfadjoint
extension $X_{\varphi _0}.$

It is known \cite{L102}, that for complex
number $\omega $
\begin{equation}
\omega =c_\infty (i)-ie^{-i\varphi _0}r_\infty (i)
\label{ad9}
\end{equation}
one has
\begin{equation}
\omega =\int_{-\infty }^\infty
\frac{\sigma _{\varphi _0}(d\lambda )}{
\lambda -i} \, . \label{ad10}
\end{equation}

The center $c_\infty (i)$ and radius $r_\infty (i)$
of the limit Weyl - Hamburger circle are given by
the formulas
\smallskip
\begin{equation}
\begin{array}{c}
c_\infty (i)=\ds{
\frac{\ds{\frac i2}-\ds{\sum_{k=0}^\infty Q_k(i;q)P_k(-i;q)}}
{\ds{\sum_{k=0}^\infty \left| P_k(i;q)\right| ^2}}}
\, ,\qquad
r_\infty (i)=
\left(
2\ds{\sum_{k=0}^\infty \left| P_k(i;q)\right|}^2
\right)^{-1}
\, .
\end{array}
\label{ad11}
\end{equation}
\smallskip
The Stieltjes transformation
\smallskip
\begin{equation}
m(z)=\int_{-\infty }^\infty
\frac{\mathrm{\sigma }_{\varphi _0}(d\lambda )}{\lambda -z}
\label{ad12}
\end{equation}
\smallskip
of the measure $\mathrm{\sigma }_{\varphi _0}$ is connected with
its value  $m(i)$ at the point $i$ by the relation \cite{L11}
\begin{equation}
m(z)=\frac{E_0(z,i)m(i)+E_1(z,i)}{D_0(z,i)m(i)+D_1(z,i)}\, .
\label{ad13}
\end{equation}
Here \,
$\ds{
E_i=\lim\limits_{n\rightarrow \infty }{E}_i^{(n)},\quad
D_i=\lim\limits_{n\rightarrow \infty }{D}_i^{(n)}\, ,
}$
\, and
\begin{equation}
\begin{array}{ccc}
{E}_0^{(n-1)}{(z,i)} & = &
    \sqrt{\left[ n\right] }\left\{
    Q_n(z;q)P_{n-1}(i;q)-P_n(i;q)Q_{n-1}(z;q)\right\}
\\[.3cm]
{E}_1^{(n-1)}{(z,i)} & = &
    \sqrt{\left[ n\right] }\left\{
    Q_n(z;q)Q_{n-1}(i;q)-Q_n(i;q)Q_{n-1}(z;q)\right\}
\\[.3cm]
{D}_0^{(n-1)}{(z,i)} & = &
    \sqrt{\left[ n\right] }\left\{
    P_n(i;q)P_{n-1}(z;q)-P_{n-1}(i;q)P_n(z;q)\right\}
\\[.3cm]
{D}_1^{(n-1)}{(z,i)} & = &
    \sqrt{\left[ n\right] }\left\{
    Q_n(i;q)P_{n-1}(z;q)-Q_{n-1}(i;q)P_n(z;q)\right\}
\end{array}
\label{i14}
\end{equation}

From given $m(z)$ one can reconstruct the spectral
measure $\mathrm{\sigma }_{\varphi _0}$ using
the inverse Stieltjes transformation \cite{a12,L11}
\begin{equation}\l{1s}
\mathrm{\sigma }_{\varphi _0}(\Delta )=
\lim\limits_{\tau \rightarrow 0}\int_\Delta
\psi (\gamma ;\tau )\mathrm{d}\gamma ;\quad
z=\gamma +i\tau ,\qquad \psi (\gamma ;\tau )=
\frac{m(z)-m(\overline{z})}{2\pi i}\, .
\end{equation}

\subsection{Harmonic oscillator}

It is wellknown that the quantum-mechanical operators
of position $X_b$ and canonically conjugate momentum
$P_b$ are selfadjoint operators with real line
as continuous spectrum. These operators are unbounded and
on the dense domain in the Hilbert space $\mathfrak{H}$
they fulfill the famous commutation relations
\begin{equation}
\left[ X_b,P_b\right] =i\hbar I,
\label{d1}
\end{equation}
ultimately connected with the Heisenberg uncertainty relation.

These operators related with creation $b^{\dagger }$ and
annihilation $b$ operators by the formulae

\begin{equation}\l{d2}
X_b=\frac 1{\sqrt{2}}\left( b^{\dagger }+b\right) ,\quad
P_b=\frac 1{i\sqrt{2}}\left( b^{\dagger }-b\right) .
\end{equation}
Operators $b^{\dagger },b$ are also unbounded, densely defined
in $\mathfrak{H}$ and adjoint to each other
$\left( (b^{\dagger })^{\dagger }\!=\!b\right).$
Consider the realization $\mathfrak{H}=\Hi$
in which the
selfadjoint number operator
$N_b=b^{\dagger }b$ $\left(
(N_b)^{\dagger }=N_b\right) $
is  diagonal. There exist the so-called vacuum state
$|0\rangle$ defined as unique state annihilated by $b,$\
$b|0\rangle =0.$
Other basis vectors are obtained by action of the creation
operator $b^{\dagger }$
\begin{equation}
|n\rangle =\frac 1{\sqrt{n!}}
\left( b^{\dagger }\right) ^n|0\rangle ,\quad
n=0,1,2,\ldots   \label{d4}
\end{equation}
On these basis states the operators
$b^{\dagger },b$ and  $N_b$ acts
according to
\begin{equation}
\begin{array}{l}
N_b\left| n\right\rangle =
   n\left| n\right\rangle ;
\\[.3cm]
b^{\dagger }\left| n\right\rangle =
   \sqrt{n+1}\left| n+1\right\rangle ,\quad
   n\geq 0;
\\[.3cm]
b\left| n\right\rangle =
    \sqrt{n}\left| n-1\right\rangle ,\quad n\geq 1,
    \qquad  b\left| 0\right\rangle =0.
\end{array}
\label{d5}
\end{equation}
This realization is known as the Fock (or number)
representation. According to the famous von Neumann
theorem this representation is unique irreducible
representation up to unitary equivalence.

The operator $X_b$ in this realization is represented
by infinite dimensional Jacobi matrix
\begin{equation}
X_b=\frac 1{\sqrt{2}}\left(
\begin{array}{cccccc}
0 & b_0 & 0 & 0 & \cdots  & \cdots  \\
b_0 & 0 & b_1 & 0 & \cdots  & \cdots  \\
0 & b_1 & 0 & b_2 & \cdots  & \cdots  \\
0 & 0 & b_2 & 0 & \ddots  & \cdots  \\
\vdots  & \vdots  & \vdots  & \ddots  & \ddots  & \ddots  \\
\vdots  & \vdots  & \vdots  & \vdots  & \ddots  & \ddots
\end{array}
\right)\, ,   \label{d6}
\end{equation}
with $b_n=\sqrt{n+1}.$ Obviously operator $X_b$
is symmetric and due to
\begin{equation}
\sum_{k=0}^\infty \frac 1{b_k}=\infty ,  \label{d7}
\end{equation}
the operator $\overline{X}_b$
is selfadjoint. It is wellknown that operators
$X_b=\frac 1{\sqrt{2}}\left(
b^{\dagger }+b\right) $ and
$P_b=\frac 1{i\sqrt{2}}\left(
b^{\dagger }-b\right) $
are related by the unitary transformation
$b^{\dagger }\ra ib.$  These operators (and their closures) are
simultaneously selfadjoint or not selfadjoint.
So, in what follows, we restrict consideration to
the case of position operator and its $q$-analogues.

In the theory of Lie groups and quantum
mechanics, special functions appear as particular matrix
elements (overlap coefficients) of appropriate operators
in corresponding representations (realizations): examples
are exponential functions, as coherent states in the
Bargmann-Fock (holomorphic) representation of \H0
\be\l{3.1} {\rm exp}({\overline{\!w}}z)=\langle w|z\rangle ,
\quad |z\rangle = \ds{{\rm e}^{z\db}}|0\rangle ,
\quad b|z\rangle =z|z\rangle ,\ee
and Hermite polynomials, as eigenvectors of the operator $N$,
in the coordinate representation,
$$H_n(x) \sim \langle n|x\rangle,
     \quad (b+\db)|x\rangle=\sqrt{2}x|x\rangle .$$

To find the generalized eigenvectors for position
operator $X_b$
\begin{equation}
X_b\left| x\right\rangle =x\left| x\right\rangle   \label{d8}
\end{equation}
we expand the state $\left| x\right\rangle $ in Fock basis
$\left| n\right\rangle $
\begin{equation}
\left| x\right\rangle =
\sum_{n=0}^\infty P_n(x)\left| n\right\rangle
\label{d9}
\end{equation}
substitute (\ref{d9}) and (\ref{d5}) into (\ref{d8}).
Then from (\ref{d5})
we obtain the following recurence relation
\begin{equation}
\sqrt{2}xP_n(x)=
\sqrt{n}P_{n-1}(x) + \sqrt{n+1}P_{n+1}(x), \qquad n\geq 0,
\label{d10}
\end{equation}
for coefficients $P_n(x)$ of (\ref{d9}) with initial conditions
\begin{equation}
P_0(x)=1,\quad P_{-1}(x)=0\, .  \label{d11}
\end{equation}
This means that
\begin{equation}
P_n(x)=\frac 1{\sqrt{2^nn!}}H_n(x),  \label{d12}
\end{equation}
where $H_n(x)$ is the usual Hermite polynomials
\be\l{d12a}
H_n(x)=(-1)^n e^{\ds{x^2}}\frac{d^n}{{dx}^n}
\left(e^{\ds{-x^2}}\right)\, ,
\ee
with recurence relation
\be\l{d13}
xH_n(x) =nH_{n-1}(x)+\frac 12H_{n+1}(x).
\ee

The basis states of the Fock representation,
i.e. eigenvectors for the number operator $N=\db b$
are represented by
\be\l{d13a2}
| n_b \rangle =
\left[ n! 2^n \sqrt{\pi} \right]^{-1/2}
H_n(x)e^{\ds{-\frac{x^2}{2}}}
\ee

The Hermite polynomials are orthogonal polynomials
\be\l{d14}
\int_{-\infty}^{\infty}\ H_m(x)\ H_n(x) d\sigma(x) =
\delta_{m n}\ {\rm d}^2_n ,
\ee
with measure
\be\l{d15}
d\sigma(x) = e^{\ds{-x^2}}\ dx
\ee
and normalization
\be\l{d16}
{\rm d}^2_n = 2^n\ n!\ \sqrt{\pi}
\ee
The classical Hamburger moment problem connected with
Jacobi matrix $X_b$ (\ref{d6}) is determined ones and its
unique solution is given by (\ref{d15}). This means that the
operator $X_b$ (and $P_b$) has zero
deficiency indices  $(0,0)$ and thus selfadjoint.

\subsection{Position operator for deformed oscillator
and $q$-Hermite polynomials}

In view of its connection with quantum groups and
algebras the deformed oscillator became the rather
popular subject in the last years (see e.g.
\cite{Z189,CK}). Here we briefly mention only main
definitions and some of the properties of
$q$-oscillator needed in the following.

Let us remark that for the  $q$-deformed oscillator
the situation became significantly richer and
interesting that in the case of the usual harmonic
oscillator recalled above. Indeed in this case
together with $q$-analogue of the Fock representation
there is plenty of inequivalent representations.
Moreover, the related Hamburger power moment problem
can be indetermine thus admitting different spectral
measures. From the point of view of the functional
analysis this means that position operator $X$ for
the $q$-oscillator may have deficiency indices $(1,1)$
and not only  $(0,0)$. In the first case the operator
$X$ has many different  selfadjoint extensions.

The deformed oscillator  algebra, \A, is generated by
three elements $a,\, \da , \, N$ with defining
relations\footnote{The relations (\ref{o1}) appears in
\cite{L3}. In \cite{L1} this relations are studied in
connection with  the generalization of the Veneziano
amplitude, by substitution of the $q$-$\Gamma$-function
instead of the standard $\Gamma$-function.}
\be\l{o1} a \da-q \da a=1\, , \ee
\be\l{o2} [N,\, a]=-a, \quad [N,\, \da ]=\da \, . \ee
Note that the generator $N$ is considered here as an
independent element, and we restrict ourselves to the
case of positive real $q\in (0,\infty)$.  The algebra
\A has a central element \cite{PK},
\be\l{o3} \zeta = q^{-N}\bigl( [N;q]-\da a \bigr)\, ,
\ee
where
\be\l{o3a}
 [N;q]:= (1-q^N)/(1-q)
\ee
is the standard basic number of $q$-analysis
(\cite{Gasper}).

In the original papers, the irreducible representation
of \A with the vacuum state $|0\rangle$
$(a|0\rangle = 0)$ was considered. The oscillator-type
representation space \H0, in the basis of eigenvectors
of the operator $N$, is
\be \l{o4}
{\cal H}_0 = \{\,|n\rangle; \quad n=0,1,2,...;
\quad a|0\rangle=0, \quad |n\rangle  =
([n;q]!)^{-1/2}(a^{\dagger})^n|0\rangle\, \}.
\ee
Due to the existence of a non-trivial central element
$\zeta$, in addition to \H0, the algebra \A has a set
of inequivalent irreducible representations $(0<q<1)$
in the spaces \Hg
$(\gamma \geq {\gamma}_c = (1-q)^{-1})$
parametrized by the value of the central element
$\zeta=-\gamma$ \cite{PK}, with the spectrum of $N$,
the set of all integers {$\Z$}.

Considering \A as an associative algebra, any
invertible transformation of the generators is
admissible; in particular, there are some natural
sets of the generators:
\cite{L6,Z189},
\be\l{o6}
A A^{\dagger}- q^{1/2}A^{\dagger}A =q^{-N/2} \, ,
\quad [N, A]=-A \, ,
\quad [N, A^{\dagger}]= A^{\dagger}, \ee
related to the quantum algebra $sl_q(2)$ via
the Schwinger realization \cite{L4,L5}, and the
following set related to the $sl_q(2)$ algebra
by a contraction procedure with fixed $q$
\cite{PK},
\be\l{o7}
[\alpha, \alpha^{\dagger}]=
q^{-N}, \quad [N,\alpha]=
-\alpha, \quad [N,\alpha^{\dagger}]
=\alpha^{\dagger} .
\ee
The equivalence of these generators is given
by the equalities
$a=q^{N/2}\alpha = q^{N/4}A$
\cite{L6,Z189}, with an obvious
one-parameter generalization, namely,
\be\l{o8}
a(\lambda)=q^{-{\frac12}\lambda N}a\, , \quad
\da(\lambda)=\da q^{-{\frac12}\lambda N}\, . \ee
This leads to the commutation relations (still one
degree of freedom)
\be\l{o9} a(\lambda) \da(\lambda) - q^{1-\lambda}
\da(\lambda) a(\lambda) = q^{-\lambda N}\, .\ee

One more formal parameter $\nu \in \R$ can be added by
a shift $N \ra N+\nu$. The corresponding set of \A
generators is denoted by $W_{p,r}^{\nu} (q)$ \cite{BDY}.
As a consequence of (\ref{o9}), namely,
\be\l{o11}
a(\lambda) \bigl( \da(\lambda)\bigr) ^m =
\bigl( p \da(\lambda) \bigr) ^m a(\lambda)  +
\bigl( p \da (\lambda) \bigr) ^{m-1} r^{N}
[m;{\frac rp}] \, ,
\ee
the normalized  basis vectors of \H0 in terms of
$ \da(\lambda)$ are given by
$$ |n\rangle = \bigl( [n;q,\lambda ]!\bigr)^{-1/2}
\bigl( \da(\lambda) \bigr) ^n |0\rangle$$
with the factorials defined as
\be \l{13} [n;q,\lambda ]!=
\prod\limits_{k=1}^{n} [k;q,\lambda ]\, ,\quad
[m;q,\lambda ] = q^{\lambda (1-m)}[m;q] \,.
\ee

The classical  moment problem  refers also to
$q$-Hermite polynomials: the latter are nothing but
polynomials of the first kind \cite{L10,L11}
for a Jacobi matrix ${\cal J}$ which is
constructed as a ``generalized coordinate'' from the
$q$-oscillator creation and annihilation operators
\cite{Z199},
\be\l{3.17}
{\cal J}(\lambda)=a(\lambda) + \da(\lambda),
\quad {\cal J}(\lambda)\, |x\rangle_{\lambda}
= 2x \,|x\rangle_{\lambda},
\ee
\be\l{3.18}
|x\rangle_{\lambda} =\sum_{n=0}^{\infty}
H_n (x;q,\lambda)|n\rangle .
\ee
Due to (\ref{3.17}), these $q$-Hermite polynomials
satisfy the following three-term recurrence relation:
\be\l{3.19} b_n(\lambda)H_{n-1}(x;q,\lambda) +
      b_{n+1}(\lambda)H_{n+1}(x;q,\lambda)=
      x\,H_n(x;q,\lambda).
\ee

The measure entering into the $q$-Hermite polynomials
$H_n (x;q,\lambda)$ orthogonality relations is connected
with the solution of the Hamburger  moment problem :
this measure is known explicitly for some cases
(see e.g.\cite{L10,L11}). This connection of the
moment problem  with Jacobi
matrices gives rise to a generalized deformation of the
oscillator identifying the matrix
$b_k \delta_{n+1,k}\, , b_k >0$
with an annihilation operator $a$. Then one gets the Wigner
commutation relation
$[a,\, \da] = F(N)$ with
$F(n) = b_{n+1}^2 - b_n^2$
and its central element
$\zeta = \bigl(c^2(N) -\da a\bigr) + const $
(see also \cite{FdP,L16,KQ}). The $q$-special
functions related to the other irreducible representations
\Hg of \A are discussed in \cite{PPK}.

Let us recall that there exists several variants
$q$-generalizations of Hermite polynomials
(see for example \cite{h1}-\cite{h7}) but in
all this definitions, based on the standard
$q$-analysis, always nonsymmetrical basic number
$[a;q]$ is used. Thus such $q$-Hermite polynomials
are little to do with symmetrical $q$-oscillator.
As far as we know our article \cite{Z199} is
the only work in which the attempt to define $q$-Hermite
polynomials related to the symmetrical $q$-oscillator is
initiated. In \cite{Z199}  we take as defining property of the
$q$-Hermite polynomials the three point recurence relation
\be\l{z0}
x H_n^q(x) = [n] H_{n-1}^q(x) + H_{n+1}^q(x) \, , \qquad
H_{0}^q(x) =1\, ,
\ee
and more physical
condition that, as in the usual case, the generating function
\be\l{z1}
g(z,x) = \sum_{n=0}^{\infty} \ds{\frac{z^n}{n!}} H_n^q(x)
\ee
is the eigenfunction of coordinate operator $X=\da+a$ in
$q$-holomorphic representation
\be\l{z2}
X=\da+a = z + {}^q\!D_z\, , \qquad  X g(z,x) = x g(z,x),
\ee
where $\qD$ is symmetrical difference $q$-derivation defined
by
\be\l{z3}
\qD f(z) =\ds{\frac{f(qz)-f(q^{-1}z)}{(q-q^{-1})z}}\, .
\ee
This gives the difference equation
\be\l{z4}
\qD g(z,x) = (x-z) g(z,x).
\ee
If we represent $g(z,x)$  in the form
\be\l{z5}
g(z,x)= 1 + g_1(z,x) + g_2(z,x) + \ldots\, ,
\ee
we obtain the sequence of equations
\be\l{z6}
\qD g_n(z,x) = (x-z) g_{n-1}(z,x)\, ,\qquad g_{0}(z,x)=1\, .
\ee
Using modification on the case of symmetric basic number of the
Jackson's $q$-integral, one easily solves  this equations.

As result one obtains
\begin{eqnarray}
g_{0}(z,x) &=& 1 \nonumber \\
g_{1}(z,x) &=& xz-\frac{\ds{z^2}}{[2]}
=\frac{\ds{xz}}{[1]!}-\frac{\ds{z^2}}{[2]!} ; \nonumber \\
g_{2}(z,x) &=& \frac{\ds{x^2z^2}}{[2]!}
-\frac{\ds{xz^3}}{[3]!} \left( [1] + [2] \right)
+\frac{\ds{z^4}}{[4]!} [3] [1] ; \nonumber \\
g_{3}(z,x) &=& \frac{\ds{x^3z^3}}{[3]!}
-\frac{\ds{x^2z^4}}{[4]!}
\left( [1] + [2] + [3] \right) + \nonumber  \\
&&+\frac{\ds{xz^5}}{[5]!}
\left( [3] [1] + [4]\left( [1] + [2]\right)\right)
- \frac{\ds{z^6}}{[6]!} [5] [3] [1] ; \l{z7}\\
g_{4}(z,x) &=& \frac{\ds{x^4z^4}}{[4]!}
-\frac{\ds{x^3z^5}}{[5]!}
\left( [1] + [2] + [3] + [4]\right) + \nonumber \\
&&+\frac{\ds{x^2z^6}}{[6]!}
\left( [3] [1] + [4]\left( [1] + [2]\right)
+ [5] \left( [1] + [2] + [3] \right) \right) - \nonumber  \\
&&- \frac{\ds{xz^7}}{[7]!} \left([5] [3] [1]
+ [6]\left([3] [1] + [4] \left([1] [2]\right)\right)\right)
+ \frac{\ds{z^6}}{[6]!} [7] [5] [3] [1] \nonumber
\end{eqnarray}
Now if we substitute (\ref{z7}) into (\ref{z5})  and join the
terms with equal degree we obtain
\begin{eqnarray}
H_{0}^q(x) &= &1 \, ,\nonumber \\
H_{1}^q(x) &= &x \, ,\nonumber \\
H_{2}^q(x) &= &x^2 -[1] \, ,\nonumber \\
H_{3}^q(x) &= &x^3 -x\left( [1] + [2] \right) \, ,\l{z8} \\
H_{4}^q(x) &= &x^4-x^2\left( [1] + [2] + [3] \right)
+[3] [1] \, ,\nonumber \\
H_{5}^q(x) &= &x^5 - x^3 \left( [1] + [2] + [3] + [4] \right)
+x \left( [3][1] + [4]\left([1] + [2]\right) \right)
\, ,\nonumber
\end{eqnarray}

It is not hard to check that
such defined $q$-Hermite polynomials
fulfil recurrent relation (\ref{z0}) which
can be used to obtain
concrete form of all other polynomials. General
expression of
$q$-Hermite polynomial looks as
\be\l{z9}
H_{n}^q(x) = x^n +
\sum_{k=1}^{\epsilon (n/2)}(-1)^k x^{n-2k} \left(
\sum_{m_k=2k-1}^{n-1} [m_k]
\sum_{m_{k-1}=2k-3}^{m_k-2} [m_{k-1}]\!
\sum_{m_{k-2}=2k-5}^{m_{k-1}-2} [m_{k-2}]
\cdot\! \ldots\! \cdot
\sum_{m_{1}=1}^{m_2-2} [m_k] \right) .
\ee
Here  $\epsilon (x)=\mathrm{Ent}(x)$ denotes the
integer part of $x.$

\subsection{The Hamburger power moment problem for
symmetricaly defor\-med oscillator}

Bellow we consider the deformed oscillator algebra \A ,
$q\in \R \;$ with fixed generators $a^{\pm },N$
fulfilling the commutation rules
\be\l{i1}
a^{-}a^{+}-qa^{+}a^{-}=q^{-N}, \qquad
\left[ N,a^{\pm }\right] =\pm a^{\pm },
\end{equation}
and usual hermiticity conditions
$\left( a^{\pm }\right) ^{\dagger }
=a^{\mp },\, N^{\dagger }=N.$

In the Fock representation action
of the generators are given by
\begin{eqnarray}
N\left| n\right\rangle &=&n\left| n
\right\rangle ;\qquad a^{+}\left|
n\right\rangle =\sqrt{\left[ n+1\right] }
\left| n+1\right\rangle ,\quad
n\geq 0;  \label{i2}
\\[5pt]
a^{-}\left| n\right\rangle &=&
\sqrt{\left[ n\right] }\left| n-1\right\rangle
,\quad n\geq 1,\qquad a^{-}\left| 0
\right\rangle =0;  \nonumber
\end{eqnarray}
where
\[
\left| n\right\rangle =
\frac 1{\sqrt{\left[ n\right] !}}\left( a^{+}\right)
^n\left| 0\right\rangle
\]
and we use the notation
\begin{equation}
\left[ \alpha \right] \equiv
\left[ \alpha \right] _q=
\frac{q^\alpha-q^{-\alpha }}{q-q^{-1}}.
\label{i3}
\end{equation}

According to these relations in the Fock space
with given above basis
$\left\{ \left| n\right\rangle \right\} _{n=0}^\infty $
the position  operator
$X:=a^{+}+a^{-}$ is described by Jacobi matrix
$X$ (\ref{d6}) in which now
$b_k=\sqrt{\left[ k+1\right] },\quad k=0,1,2,\ldots .$

From eigenvalue equation
$X\left| x\right\rangle =x\left| x\right\rangle\, ,$
$\left| x\right\rangle
=\sum_{n=0}^\infty P_n(x;q)\left| n\right\rangle\, ,$
one obtains following recurrent relations for
q-Hermitian polynomials $P_n(x;q)$
\begin{equation}
\begin{array}{c}
\sqrt{\left[ n\right] }P_{n-1}(x;q)+\sqrt{\left[ n+1\right] }
P_{n+1}(x;q)=xP_n(x;q),\quad n\geq 1 \\[5pt]
P_1(x;q)=xP_0(x;q),\qquad P_0(x;q)=1.
\end{array}
\label{i6}
\end{equation}
These q-Hermitian polynomials $P_n(x;q)$ are polynomials
of the 1-st kind for Jacobi matrix $X.$
Note that polynomials $Q_n(x;q)$ of the
2-nd kind for the same Jacobi matrix fulfil
the same recurrence relations (\ref{i6}) but with
other initial conditions
\[
Q_1(x;q)=1,\qquad Q_0(x;q)=0.
\]
The wellknown connection \cite{L10}\
between polynomials of the 1-st and 2-nd
kind looks in our case as
\begin{equation}
P_{n-1}(x;q)Q_n(x;q)-P_n(x;q)Q_{n-1}(x;q)=
\ds{\frac 1{\sqrt{\left[ n\right]}}}
\label{i7}
\end{equation}

Because for $q>0,q\neq 1$ we have
\[
\sqrt{\left[ n\right]} \sqrt{\left[ n+2\right] }
\leq \left[ n+1\right] ,\quad n\geq 0 ,\qquad
\sum_{n=0}^\infty \left[ n+1\right] ^{-\frac 12}<\infty
\]
then, according to (\cite{L11}, thm 1.5)
the closure $\overline{X}$ of the position operator
$X\;$ is only closed but not selfadjoint operator and
has deficiency indices $\left( 1.1\right) .$ This
means that coordinate operator $X$ has family of
selfadjoint extensions.

The related Hamburger moment problem is
indetermine one, and there exist family of
measures $\sigma $ fulfilling relations (\ref{ad7}).
The spectral measure $\sigma _{\varphi _0}$ gives
the ''extremal'' solution of the Hamburger moment
problem (\ref{ad7}) related
to the selfadjoint extension $X_{\varphi _0}.$

{\bf Remark}.  The Stieltjes transformation $m(z)$
of the spectral measure $\sigma _{\varphi_0 }$
one can obtain by the relation \cite{L10}
\begin{equation}
m(z)=\frac{{A}{(z)t-C}{(z)}}{{B}{(z)t-D}{(z)}}  \label{i15}
\end{equation}
where the elements of Nevanlinna matrix
\begin{equation}
\left(
\begin{array}{cc}
A(z) & C(z) \\
B(z) & D(z)\\
\end{array}
\right)
\end{equation}
are connected with ${E}_i^{(n)}$ and ${D}_i^{(n)}$
by the  relations
\begin{equation}
\begin{array}{ccccc}
{A}{(z)}=&{E}_1{(z;0)}&\quad&{B}{(z)}=&-{D}_1{(z;0)}
\\
{C}{(z)}=&{E}_0{(z;0)}&&{D}{(z)}=&-{D}_0{(z;0)}
\end{array}
\label{i16}
\end{equation}
Unfortunately, there is no way of finding at once
the value $t$ in (\ref{i15}) which  correspond
to the considered selfadjoint extension of $X$,
labelled by $\varphi_0$. To obtain
the Stieltjes transformation $m(z)$  of
the spectral measure $\sigma _{\varphi_0 }$
we need to use the more complicated relation (\ref{ad13}).

As it was mentioned above the inverse Stieltjes
transformation \cite{a12,L11} allows to
reconstruct the spectral measure
$\mathrm{\sigma }_{\varphi _0}$
from given $m(z)$ according to formula (\ref{1s}).

\section{Stieltjes transformation $m(z)$ of the
spectral measure}
\setcounter{equation}{0}

\subsection{Computation of $m(i)$}

In this section we will get the explicit form of the
Stieltjes transformation $m(i)$  of the spectral measure
$\sigma _{\varphi_0 }$
in term of the elements of Jacobi matrix $X.$
In this connection we need to consider some non-standard
$q$-series, which are possibly of a special interest by
themselves. The polynomials $P_n$ and $Q_n$ can be
represented in the form
\begin{equation}
P_n(x;q)=\sum_{m=0}^{\epsilon (\frac n2)}
\frac{(-1)^m}{\sqrt{\left[ n\right]!}}
\alpha _{2m-1,n-1}x^{n-2m}  \label{a1.1}
\end{equation}
\begin{equation}\l{a1.2}
\alpha _{\,-1;n-1}\equiv 1;\qquad \alpha
_{2m-1;n-1}=\sum\limits_{k_1=2m-1}^{n-1}\left[ k_1\right]
\sum\limits_{k_2=2m-3}^{k_1-2}\left[ k_2\right] \ldots
\sum\limits_{k_m=1}^{k_{m-1}-2}\left[ k_m\right] ,\quad m\geq 1
\end{equation}
\begin{equation}
Q_{n+1}(x;q)=\sum_{m=0}^{\epsilon (\frac n2)}
\frac{(-1)^m}{\sqrt{\left[n+1\right] !}}
\beta _{2m,n}x^{n-2m}  \label{a1.3}
\end{equation}
\begin{equation}
\beta _{0;n}\equiv 1;\qquad \beta _{2m;n}=
\sum\limits_{k_1=2m-1}^n\left[k_1\right]
\sum\limits_{k_2=2m-2}^{k_1-2}\left[ k_2\right]
\ldots \sum\limits_{k_m=2}^{k_{m-1}-2}
\left[ k_m\right] ,\quad m\geq 1 , \label{a1.4}
\end{equation}
recall that $\epsilon (x)=\mathrm{Ent}(x)$ denotes the
integer part of $x.$ We remark that the relations
(\ref{a1.1})-(\ref{a1.4}) give us the non-standard
representation for the polynomials $P_{n}$ and
$Q_{n}$. So we need to obtain further some properties
of the coefficients $\alpha _{m;n}$ and $\beta _{m;n}$ .

From recurrence relations for $P_n$ and $Q_n$
(\ref{i6}) it follows that the coefficients
$\alpha _{m;n}$ and $\beta _{m;n}$ fulfills the relations
\begin{equation}
\alpha _{2m-1;n}=\left[ n\right]
\alpha _{2m-3;n-2}+\alpha _{2m-1;n-1};
\label{a1.5}
\end{equation}
\begin{equation}
\beta _{2m;n}=\left[ n\right]
\beta _{2m-2;n-2}+\beta _{2m;n-1}.
\label{a1.6}
\end{equation}

From the relations (\ref{ad9}), (\ref{ad10}) and
(\ref{ad12}) we have
\begin{equation}
m(i)=\int_{-\infty }^\infty
\frac{\mathrm{\sigma }_{\varphi _0}(d\lambda )}
{ \lambda -i}=c_\infty (i)-ie^{-i\varphi _0}r_\infty (i).
\label{a1.7}
\end{equation}
Thus to compute $m(i)$ we must compute the center
$c_\infty (i)$ and the  radius $r_\infty (i)$ of
Weyl - Hamburger circle at a point $z=i.$

Using the formula (5) from
(\cite{L10}, ch.1, \S 2) and (\ref{a1.1}) we have
\begin{equation}
\begin{array}{cl}
\sum_{k=0}^{n-1}\left| P_k(i;q)\right| ^2= &
\frac{\sqrt{\left[ n\right] }}{
2i} \left( P_n(i;q)P_{n-1}(\overline{i};q)-
P_{n-1}(i;q)P_n(\overline{i};q)\right)
\\[8pt]
& =\ds{\frac 1{\left[ n-1\right] !}}
\left( \sum_{m=0}^{\epsilon (\frac n2)}
\alpha _{2m-1,n-1}\right)
\left( \sum_{m=0}^{\epsilon (\frac{n-1}2)}
\alpha _{2m-1,n-2}\right)
\end{array}
\label{a1.8}
\end{equation}
Let us introduce the auxiliary functions
\begin{equation}
\Psi (q)=\lim\limits_{s\rightarrow \infty }
\Psi _s(q),\qquad\quad \Psi _s(q)=1+\sum_{k=1}^s
\frac{\left[ 2k-1\right] !!}{\left[ 2k\right] !!}
\label{a1.9a}
\end{equation}
\begin{equation}
\Phi (q)=\lim\limits_{s\rightarrow \infty }
\Phi _s(q),\qquad\quad \Phi_s(q)=1+\sum_{k=1}^s
\frac{\left[ 2k\right] !!}{\left[ 2k+1\right] !!}
\label{a1.9b}
\end{equation}
Here
\[
\left[ 2n\right] !!=\left[ 2n\right] \left[ 2n-2\right]
\cdot \ldots \cdot \left[ 2\right] ,\qquad
\left[ 2n-1\right] !!=\left[ 2n-1\right]
\left[2n-3\right] \cdot \ldots \cdot \left[ 1\right] .
\]
It is not difficult to check that for $q>0,$ q$\neq 1$
the functions $\Psi(q) $ (\ref{a1.9a}) and $\Phi (q)$
(\ref{a1.9b}) are well defined as convergent q-series.
In fact it is follow from the inequality
$[2s]\geq ([2s-1][2s+1])^{\frac 12}$ that
\begin{equation}
{\frac{[2k-1]!!}{[2k]!!}}\leq {\frac{[1]^{\frac{1}{2}}}
{[2k+1]^{\frac{1}{2}}}}
\end{equation}
In view of $\sum_{n=0}^\infty
\left[ n+1\right] ^{-\frac 12}<\infty$
for $q>0$, $q\neq 1$ we have that $\Psi (q)$ is well
defined. The same is true of $\Phi (q)$.

It is convenient to take out the leading terms
$\alpha _{2\epsilon (\frac n2)-1,n-1}$ and
$\alpha _{2\epsilon (\frac{n-1}2)-1,n-2}$ from the
sums in the first and the second parenthesises in the
right hand side of the expression (\ref{a1.8}). This gives
us the common factor equal to
\begin{equation}
\begin{array}{cllcl}
1) & \alpha _{2p-1,2p}\alpha _{2p-1,2p-1}&=
\left[ 2p\right] !\Psi _p(q), & \mathrm{if} & n=2p+1 \\
2) & \alpha _{2p-1,2p-1}\alpha _{2p-3,2p-2}&=
\left[ 2p-1\right] !\Psi _{p-1}(q), & \mathrm{if} & n=2p
\end{array}
\label{a1.10}
\end{equation}
where we take into account that from (\ref{a1.2}) and
$\Phi (q)$ (\ref{a1.9b}) we have
\begin{equation}
\alpha _{2p-1,2p}=\left[ 2p\right] !!\,
\Psi _p(q),\qquad \alpha _{2p-1,2p-1}=
\left[ 2p-1\right] !!\,.  \label{a1.11}
\end{equation}
This allows us rewrite (\ref{a1.8}) in the form
\begin{equation}
\begin{array}{cc}
1) & \sum\limits_{k=0}^{2p}\left| P_k(i;q)
\right| ^2 = \Psi_p(q)\sum\limits_{m=0}^p
\ds{\frac{\alpha _{2m-1,2p}}{\alpha _{2p-1,2p}}}
\sum\limits_{m=0}^p
\ds{\frac{\alpha _{2m-1,2p-1}}{\alpha _{2p-1,2p-1}}}
\\[5pt]
2) & \sum\limits_{k=0}^{2p-1}\left| P_k(i;q)
\right| ^2=\Psi _{p-1}(q)\sum\limits_{m=0}^p
\ds{\frac{\alpha _{2m-1,2p-1}}{\alpha _{2p-1,2p-1}}}
\sum\limits_{m=0}^{p-1}
\ds{\frac{\alpha _{2p-1,2p-2}}{\alpha_{2p-3,2p-2}}}
\end{array}
\label{a1.12}
\end{equation}
With help of the recurrent relations (\ref{a1.5}) one obtain
\begin{equation}
\begin{array}{c}
(A^{(1)}_\alpha )_{p}:=\sum\limits_{m=0}^p
\ds{\frac{\alpha _{2m-1,2p-1}}{\alpha _{2p-1,2p-1}}}
=1+\sum\limits_{k=1}^p
\ds{\frac 1{\left[ 2k-1\right]!!}}
\sum\limits_{m=0}^{k-1}\alpha _{2m-1,2k-2}.
\\[5pt]
(A^{(2)}_\alpha)_{p}:=
\sum\limits_{m=0}^p\ds{
\frac{\alpha _{2m-1,2p}}{\alpha _{2p-1,2p}}}=
1+\ds{\frac 1{\Psi _p(q)}}
\sum\limits_{k=1}^p
\ds{\frac 1{\left[ 2k\right] !!}}
\sum\limits_{m=0}^{k-1}\alpha _{2m-1,2k-1};
\end{array}
\label{a1.13}
\end{equation}
It is clear, that $(A^{(i)}_\alpha)_{p}$, $i=1,2$,
are the positive increasing quantities. So exist the
limits $A^{(i)}_\alpha :=\lim\limits_{p\rightarrow \infty}
(A^{(i)}_\alpha)_{p}$, $i=1,2$. It follows from
(\ref{a1.12}) and the inequalities
\begin{equation}
\sum\limits_{k=0}^{\infty}\left| P_k(i;q)
\right| ^2<\infty,\qquad \Psi (q)<\infty,\quad(q>0,q\neq 1),
\end{equation}
that the expressions
\begin{equation}
A^{(1)}_\alpha =1+\sum\limits_{k=1}^\infty
\frac 1{\left[ 2k-1\right]!!}\sum\limits_{m=0}^{k-1}
\alpha _{2m-1,2k-2},  \label{a1.14}
\end{equation}
\begin{equation}
A_\alpha ^{(2)}=1+\frac 1{\Psi (q)}\sum\limits_{k=1}^\infty
\frac 1{\left[2k\right] !!}\sum\limits_{m=0}^{k-1}
\alpha _{2m-1,2k-1},  \label{a1.15}
\end{equation}
are well defined for $q>0,$ $q\neq 1.$

Finally from the relations
(\ref{ad11}), (\ref{a1.12}) - (\ref{a1.15})
we obtain the following expression for the radius
of limit Weyl - Hamburger circle
\begin{equation}
r_\infty (i)=
 \left( {2 A_\alpha ^{(1)}\Psi (q)A_\alpha ^{(2)} }
\right) ^{-1} \, .
\label{a1.16}
\end{equation}

Now we go to consideration of the center $c_\infty
(i)=\lim\limits_{n\rightarrow \infty }c_n(i),$
 where (see \cite{L10})
\begin{equation}
c_n(i)=-\frac{D_1^{(n-1)}(-i;i)}{2i\;
\sum\limits_{k=0}^{n-1}\left|
P_k(i)\right| ^2}  \label{a1.17}
\end{equation}
For terms of the numerator of this expression
\begin{equation}
D_1^{(n-1)}(-i;i)=A(n)+B(n)  \label{a1.18}
\end{equation}
in view of the relations (\ref{i14}), (\ref{a1.1}) and
(\ref{a1.3}) one obtain
\[
A(n)=-\sqrt{\left[ n\right] }Q_{n-1}(i)P_n(-i)=
\frac 1{\left[ n-1\right]
!}\left( \sum_{m=0}^{\epsilon (\frac{n-2}2)}
\beta _{2m,n-2}\right) \left(
\sum_{m=0}^{\epsilon (\frac n2)}
\alpha _{2m-1,n-1}\right)
\]
\[
B(n)=\sqrt{\left[ n\right] }Q_n(i)P_{n-1}
(\overline{i})=\frac 1{\left[
n-1\right] !}\left( \sum_{m=0}^{\epsilon
(\frac{n-1}2)}\beta _{2m,n-1}\right) \left(
\sum_{m=0}^{\epsilon (\frac{n-1}2)}\alpha
_{2m-1,n-2}\right)
\]

From (\ref{a1.4}) and (\ref{a1.9b}) one obtain also
\begin{equation}
\beta _{2p,2p+1}=\left[ 2p+1\right] !!\Phi _p(q),\qquad \beta
_{2p,2p}=\left[ 2p\right] !!  \label{a1.19}
\end{equation}

Using (\ref{a1.11}) and (\ref{a1.19}) and repeat
the calculations used above in consideration of the
expression (\ref{a1.16}) we obtain
\begin{equation}
\begin{array}{c}
\lim\limits_{p\rightarrow \infty }A(2p)=
\lim\limits_{p\rightarrow \infty}B(2p+1)=
A_\alpha ^{(1)}A_\beta ^{(1)}; \\[3pt]
\lim\limits_{p\rightarrow \infty }A(2p+1)=
\lim\limits_{p\rightarrow \infty
}B(2p)=\Psi (q)\Phi (q)A_\alpha ^{(2)}A_\beta ^{(2)},
\end{array}
\label{a1.20}
\end{equation}
where together with (\ref{a1.14}) and (\ref{a1.15}) we
use the notation
\begin{equation}
A_\beta ^{(1)}=1+\sum\limits_{k=1}^\infty
\frac 1{\left[ 2k\right]!!}\sum_{m=0}^{k-1}\beta _{2m,2k-1}
\label{a1.21}
\end{equation}
\begin{equation}
A_\beta ^{(2)}=1+\frac 1{\Phi (q)}\sum\limits_{k=2}^\infty
\frac 1{\left[ 2k-1\right] !!}\sum_{m=0}^{k-2}
\beta _{2m,2k-2}  \label{a1.22}
\end{equation}
It is not difficult to check that $A_\beta ^{(1)}$
and $A_\beta ^{(2)}$ are well defined for $q>0,$
 $q\neq 1.$

Finally from the relations (\ref{a1.16})-(\ref{a1.18}) and
(\ref{a1.20}) we obtain
\begin{equation}
C_\infty (i)=\frac i2\frac{\left( A_\alpha ^{(1)}
A_\beta ^{(1)}+\Psi(q)A_\alpha ^{(2)}\Phi (q)
A_\beta ^{(2)}\right) }{A_\alpha ^{(1)}\Psi(q)
A_\alpha ^{(2)}}  \label{a1.23}
\end{equation}
From (\ref{a1.7}), (\ref{a1.16}) and (\ref{a1.23}) we have
\begin{equation}
m(i)=\frac{-\sin \varphi _0}{A_\alpha ^{(1)}
\Psi (q)A_\alpha ^{(2)}}+
i\frac{-\cos \varphi _0+A_\alpha ^{(1)}A_\beta ^{(2)}+
\Psi (q)A_\alpha ^{(2)}\Phi(q)A_\beta ^{(2)}}
{A_\alpha ^{(1)}\Psi (q)A_\alpha ^{(2)}}  \label{a1.24}
\end{equation}

\textbf{Remark 1.}Note that from (\ref{a1.1}) and
(\ref{a1.3}) it follows that
\begin{equation}
\begin{array}{c}
Q_{n-1}(i;q)P_n(i;q)=(-1)^nQ_{n-1}(i;q)P_n(-i;q)
=(-1)^{n-1}A(n)\left[ n\right] ^{-\frac 12} \\[5pt]
Q_n(i;q)P_{n-1}(i;q)=(-1)^{n-1}Q_n(i;q)P_{n-1}(-i;q)
=(-1)^{n-1}B(n)\left[ n\right] ^{-\frac 12}
\end{array}
\label{a1.25}
\end{equation}

Then from (\ref{i7}), (\ref{a1.20}) and (\ref{a1.25})
we have
\begin{equation}
W=-\Psi (q)A_\alpha ^{(2)}\Phi (q)A_\beta ^{(2)}+
A_\alpha ^{(1)}A_\beta ^{(1)}=1
\label{a1.26}
\end{equation}
\textbf{Remark 2.} We emphasize that all results of
the present section are hold for every Jacobi matrix
(\ref{d6}) with restrictions
\[
b_ {n} b_{n+2}\ \leq \left(b_{n+1}\right)^2,\quad n\geq 0;
\qquad
\sum_{n=0}^\infty \left(b_{n}\right)^{-\frac 12}<\infty
\]
which provide the associated moment problem to be indetermine.

\subsection{Computation of the $m(z)$}

For computation of the $m(z)$ (\ref{ad13}) we must find
the functions (\ref {i14}). To this end we represent
the complex variable $z$ in the polar form
\begin{equation}
s=q^\kappa e^{i\varphi },\quad q>0,\;q\neq 1,\quad
\kappa \in \R \quad , \quad 0\leq \varphi \leq 2\pi
\label{a2.1}
\end{equation}

For computation the functions (\ref{i14}) we use the
same reasoning as in the sect.1 for obtaining the
formula (\ref{a1.16}). For doing so we must
introduce together with auxiliary quantities
$A_\varepsilon ^{(i)}$ ($
\varepsilon =\alpha ,\beta ;$ $i=1,2)$
also the following expressions
\begin{equation}
A_{\alpha ,\kappa }^{(1)}(\varphi ):=
1+\sum\limits_{k=1}^\infty \frac1{\left[ 2k-1\right] !!}
\sum\limits_{m=0}^{k-1}(-1)^{k-m}
\left( q^\kappa e^{i\varphi }\right) ^{2(k-m)}
\alpha _{2m-1,2k-2};  \label{a2.2}
\end{equation}
\begin{equation}
A_{\alpha ,\kappa }^{(2)}(\varphi ):=
1+\frac 1{\Psi (q)}\sum\limits_{k=1}^\infty
\frac 1{\left[ 2k\right]!!}
\sum\limits_{m=0}^{k-1}(-1)^{k-m}
\left( q^\kappa e^{i\varphi }\right)
^{2(k-m)}\alpha _{2m-1,2k-1};  \label{a2.3}
\end{equation}
\begin{equation}
A_{\beta ,\kappa }^{(1)}(\varphi ):=
1+\sum\limits_{k=1}^\infty \frac 1{\left[ 2k\right] !!}
\sum\limits_{m=0}^{k-1}(-1)^{k-m}
\left( q^\kappa e^{i\varphi }\right) ^{2(k-m)}
\beta _{2m,2k-1}  \label{a2.4}
\end{equation}
\begin{equation}
A_{\beta ,\kappa }^{(2)}(\varphi ):=
1+\frac 1{\Phi (q)}\sum\limits_{k=2}^\infty
\frac 1{\left[ 2k-1\right]!!}
\sum\limits_{m=0}^{k-2}(-1)^{k-m-1}
\left( q^\kappa e^{i\varphi }\right)
^{2(k-m-1)}\beta _{2m,2k-2};  \label{a2.5}
\end{equation}
Note that $A_{\varepsilon ,\kappa =0}^{(i)}(\varphi
=\frac \pi 2)=A_\varepsilon ^{(i)}.$
In terms of the above notations we have
\begin{equation}
E_0(z;i)=ie^{i\varphi }q^\kappa \Psi (q)
A_\alpha ^{(2)}\Phi (q)A_{\beta,\kappa }^{(2)}
(\varphi )+A_\alpha ^{(1)}
A_{\beta ,\kappa }^{(1)}(\varphi );
\label{a2.6}
\end{equation}
\begin{equation}
E_1(z;i)=-iA_{\beta ,\kappa }^{(1)}(\varphi )
\Phi (q)A_\beta^{(2)}+e^{i\varphi }q^\kappa
A_\beta ^{(1)}\Phi (q)A_{\beta ,\kappa }^{(2)}(\varphi );
\label{a2.7}
\end{equation}
\begin{equation}
D_1(z;i)=ie^{i\varphi }q^\kappa \Psi (q)
A_{\alpha ,\kappa }^{(2)}(\varphi )\Phi (q)A_\beta ^{(2)}+
A_\beta ^{(1)}A_{\alpha ,\kappa }^{(1)}(\varphi );
\label{a2.8}
\end{equation}
\begin{equation}
D_0(z;i)=i\Psi (q)A_\alpha ^{(2)}A_{\alpha ,\kappa }^{(1)}
(\varphi )-e^{i\varphi }q^\kappa A_\alpha ^{(1)}\Psi (q)
A_{\alpha ,\kappa }^{(2)}(\varphi );
\label{a2.9}
\end{equation}

Let us denote
\begin{equation}
A_{\varepsilon ,\kappa }^{(i)}(\varphi )
=A_{\varepsilon ,\kappa }^{(i)R}(\varphi )+
iA_{\varepsilon ,\kappa }^{(i)I}(\varphi )
\label{a2.10}
\end{equation}
Using the relations (\ref{a1.26}) and (\ref{a2.10})
from (\ref{i15}), (\ref{a1.24}) and (\ref{a2.6}) - (\ref{a2.9})
one obtains
\[
m(z)=\frac{\mathcal{N}_{m(z)}}{\mathcal{D}_{m(z)}}
\]
\begin{equation}
\begin{array}{cc}
\mathcal{N}_{m(z)}= & \left\{
-\ds{\frac 12}\mathcal{N}_1-
\ds{\frac{\sin \varphi _0}2} \mathcal{N}_2+
\ds{\frac{\cos \varphi _0}2} \mathcal{N}_3\right\} -
\\[8pt]
& i\left\{ -\frac 12\mathcal{N}_4+
\ds{\frac{\sin \varphi _0}2} \mathcal{N}_5+
\ds{\frac{\cos \varphi _0}2} \mathcal{N}_6\right\} ;
\end{array}
\label{a2.11}
\end{equation}
\medskip
\begin{equation}
\begin{array}{cc}
\mathcal{D}_{m(z)}= & \left\{
-\ds{\frac 12}\mathcal{D}_1-
\ds{\frac{\sin \varphi _0}2} \mathcal{D}_2+
\ds{\frac{\cos \varphi _0}2} \mathcal{D}_3\right\} -
\\[8pt]
& i\left\{ -\ds{\frac 12}\mathcal{D}_4+
\ds{\frac{\sin \varphi _0}2} \mathcal{D}_5+
\ds{\frac{\cos \varphi _0}2} \mathcal{D}_6\right\} ,
\end{array}
\label{a2.12}
\end{equation}
where
\[
\mathcal{N}_1=\frac{A_{\beta ,\kappa }^{(1)I}(\varphi )}
{\Psi (q)A_\alpha ^{(2)}}-q^\kappa \cos \varphi
\frac{\Phi (q)A_{\beta ,\kappa }^{(2)R}(\varphi )}
{A_\alpha ^{(1)}}+q^\kappa \sin \varphi
\frac{\Phi (q)A_{\beta ,\kappa }^{(2)I}(\varphi )}
{A_\alpha ^{(1)}}
\]
\[
\mathcal{N}_2=\frac{A_{\beta ,\kappa }^{(1)R}(\varphi )}
{\Psi (q)A_\alpha ^{(2)}}-q^\kappa \cos \varphi
\frac{\Phi (q)A_{\beta ,\kappa}^{(2)I}(\varphi )}
{A_\alpha ^{(1)}}-q^\kappa \sin \varphi \frac{\Phi
(q)A_{\beta ,\kappa }^{(2)R}(\varphi )}{A_\alpha ^{(1)}}
\]
\[
\mathcal{N}_3=\frac{A_{\beta ,\kappa }^{(1)I}(\varphi )}
{\Psi (q)A_\alpha ^{(2)}}-q^\kappa \cos \varphi
\frac{\Phi (q)A_{\beta ,\kappa }^{(2)R}(\varphi )}
{A_\alpha ^{(1)}}+q^\kappa \sin \varphi
\frac{\Phi (q)A_{\beta ,\kappa }^{(2)I}(\varphi )}
{A_\alpha ^{(1)}}
\]
\[
\mathcal{N}_4=\frac{A_{\beta ,\kappa }^{(1)R}(\varphi )}
{\Psi (q)A_\alpha ^{(2)}}+q^\kappa \cos \varphi
\frac{\Phi (q)A_{\beta ,\kappa }^{(2)I}(\varphi )}
{A_\alpha ^{(1)}}+q^\kappa \sin \varphi \frac{\Phi (q)
A_{\beta ,\kappa }^{(2)R}(\varphi )}{A_\alpha ^{(1)}}
\]
\[
\mathcal{N}_5=\frac{A_{\beta ,\kappa }^{(1)I}(\varphi )}
{\Psi (q)A_\alpha ^{(2)}}-q^\kappa \cos \varphi
\frac{\Phi (q)A_{\beta ,\kappa }^{(2)I}(\varphi )}
{A_\alpha ^{(1)}}+q^\kappa \sin \varphi \frac{\Phi (q)
A_{\beta ,\kappa }^{(2)R}(\varphi )}{A_\alpha ^{(1)}}
\]
\[
\mathcal{N}_6=\frac{A_{\beta ,\kappa }^{(1)R}(\varphi )}
{\Psi (q)A_\alpha ^{(2)}}-q^\kappa \cos \varphi
\frac{\Phi (q)A_{\beta ,\kappa }^{(2)I}(\varphi )}
{A_\alpha ^{(1)}}-q^\kappa \sin \varphi \frac{\Phi (q)
A_{\beta ,\kappa }^{(2)R}(\varphi )}
{A_\alpha ^{(1)}}
\]
\[
\mathcal{D}_1=\frac{A_{\alpha ,\kappa }^{(1)R}(\varphi )}
{A_\alpha ^{(1)}}+q^\kappa \cos \varphi
\frac{\Psi (q)A_{\alpha ,\kappa }^{(2)I}(\varphi )}
{\Psi (q)A_\alpha ^{(2)}}+q^\kappa \sin \varphi
\frac{\Psi (q)A_{\alpha ,\kappa }^{(2)R}(\varphi )}
{\Psi (q)A_\alpha ^{(2)}}
\]
\[
\mathcal{D}_2=\frac{A_{\alpha ,\kappa }^{(1)I}(\varphi )}
{A_\alpha ^{(1)}}-q^\kappa \cos \varphi
\frac{\Psi (q)A_{\alpha ,\kappa }^{(2)R}(\varphi )}
{\Psi (q)A_\alpha ^{(2)}}+q^\kappa \sin \varphi
\frac{\Psi (q)A_{\alpha ,\kappa }^{(2)I}(\varphi )}
{\Psi (q)A_\alpha ^{(2)}}
\]
\[
\mathcal{D}_3=\frac{A_{\alpha ,\kappa }^{(1)R}
(\varphi )}{A_\alpha ^{(1)}}
-q^\kappa \cos \varphi
\frac{\Psi (q)A_{\alpha ,\kappa }^{(2)I}(\varphi )}{
\Psi (q)A_\alpha ^{(2)}}-q^\kappa
\sin \varphi \frac{\Psi (q)A_{\alpha
,\kappa }^{(2)R}(\varphi )}{\Psi (q)A_\alpha ^{(2)}}
\]
\[
\mathcal{D}_4=\frac{A_{\alpha ,\kappa }^{(1)I}
(\varphi )}{A_\alpha ^{(1)}}
-q^\kappa \cos \varphi
\frac{\Psi (q)A_{\alpha ,\kappa }^{(2)R}(\varphi )}{
\Psi (q)A_\alpha ^{(2)}}+
q^\kappa \sin \varphi \frac{\Psi (q)A_{\alpha
,\kappa }^{(2)I}(\varphi )}{\Psi (q)A_\alpha ^{(2)}}
\]
\[
\mathcal{D}_5=\frac{A_{\alpha ,\kappa }^{(1)R}
(\varphi )}{A_\alpha ^{(1)}}
-q^\kappa \cos \varphi \frac{\Psi (q)
A_{\alpha ,\kappa }^{(2)I}(\varphi )}{
\Psi (q)A_\alpha ^{(2)}}-q^\kappa \sin
\varphi \frac{\Psi (q)A_{\alpha
,\kappa }^{(2)R}(\varphi )}{\Psi (q)A_\alpha ^{(2)}}
\]
\[
\mathcal{D}_6=-\frac{A_{\alpha ,\kappa }^{(1)I}
(\varphi )}{A_\alpha ^{(1)}}
-q^\kappa \cos \varphi \frac{\Psi (q)
A_{\alpha ,\kappa }^{(2)R}(\varphi )}{
\Psi (q)A_\alpha ^{(2)}}+
q^\kappa \sin \varphi \frac{\Psi (q)A_{\alpha
,\kappa }^{(2)I}(\varphi )}{\Psi (q)A_\alpha ^{(2)}}
\]

The formulas analogous to (\ref{a2.11}) and
(\ref{a2.12}) for $m(\overline{z})$ can be obtained
from (\ref{a2.11}) and (\ref{a2.12}) by replacement $
\varphi \rightarrow -\varphi $ with help of the relations
\begin{equation}
A_{\varepsilon ,\kappa }^{(i)I}(-\varphi )
=-A_{\varepsilon ,\kappa
}^{(i)I}(\varphi ),\quad i=1,2;\quad
\varepsilon =\alpha ,\beta .
\label{a2.13}
\end{equation}
\section{Support of the spectral measure
$\sigma _{\varphi _0}$}
\setcounter{equation}{0}

From (\ref{a2.11}) and (\ref{a2.12}) and their analogues for
$m(\overline{z})$ we check that Re($m(z)-m(\overline{z})$)=0.
Thus for reconstruction of the spectral measure
$\sigma _{\varphi _0}$ from its Stieltjes transformation
$ m(z)$ according to the relation (\ref{i15}) we must find
Im($m(z)-m( \overline{z})$). We have
\begin{equation}
\begin{array}{c}
\mathrm{Im}(m(z)-m(\overline{z}))=
\qquad \qquad \qquad \qquad \qquad \qquad \\[10pt]
\qquad \qquad
2\ds{\frac{q^\kappa \left( \mathrm{Im}
\mathcal{B}_{\alpha,\beta,\kappa }^{(+)}+
\cos \varphi _0\mathrm{Im}
\mathcal{B}_{\alpha ,\beta ,\kappa }^{(-)}\right)+
\sin \varphi _0\mathrm{Im}
\mathcal{A}_{\alpha ,\beta ,\kappa }}
{\mathcal{A}_{\alpha ,\kappa }^{(+)}+
\cos \varphi _0\mathcal{A}_{\alpha ,\kappa }^{(-)}+
2q^\kappa \sin \varphi _0\mathrm{Re}\mathcal{B}_\alpha }}\, ,
\end{array}
\label{a3.5}
\end{equation}
where we use the following notation
\begin{equation}
\mathcal{A}_{\alpha ,\kappa }^{(\pm )}:=
\left( \frac{A_{\alpha ,\kappa
}^{(1)}(\varphi )}{A_\alpha ^{(1)}}\right) ^2\pm
q^{2\kappa }\left( \frac{
\Psi (q)A_{\alpha ,\kappa }^{(2)}(\varphi )}
{\Psi (q)A_\alpha ^{(2)}}\right)
^2  \label{a3.1}
\end{equation}
\begin{equation}
\mathcal{A}_{\alpha ,\beta ,\kappa }:=\frac{A_{\alpha ,\kappa
}^{(1)}(\varphi )\overline{A_{\beta ,\kappa }^{(1)}
(\varphi )}+q^{2\kappa }
\overline{\Psi (q)A_{\beta ,\kappa }^{(2)}
(\varphi )}\Phi (q)A_{\beta
,\kappa }^{(2)}(\varphi )}
{A_\alpha ^{(1)}\Psi (q)A_\alpha ^{(2)}}
\label{a3.2}
\end{equation}
\begin{equation}
\mathcal{B}_{\alpha ,\beta ,\kappa }^{(\pm )}:
=e^{i\varphi }\left( \frac{
\overline{A_{\alpha ,\kappa }^{(1)}
(\varphi )}\Phi (q)A_{\beta ,\kappa}^{(2)}(\varphi )}
{\left( A_\alpha ^{(1)}\right) ^2}\pm \frac{\overline{
A_{\beta ,\kappa }^{(1)}(\varphi )}\Psi (q)
A_{\alpha ,\kappa }^{(2)}(\varphi
)}{\left( \Psi (q)A_\alpha ^{(2)}\right) ^2}\right)
\label{a3.3}
\end{equation}
\begin{equation}
\mathcal{B}_\alpha :=
e^{i\varphi }\frac{\overline{A_{\alpha ,\kappa
}^{(1)}(\varphi )}
\Psi (q)A_{\alpha ,\kappa }^{(2)}(\varphi )}{A_\alpha
^{(1)}\Psi (q)A_\alpha ^{(2)}}  \label{a3.4}
\end{equation}
Recall that according to (\ref{i15})
\begin{equation}
\psi (\gamma ;\tau )=\frac{m(z)-m(\overline{z})}
{2\pi i}=\frac{\mathcal{N}
_\psi }{\mathcal{D}_\psi }  \label{a3.5a}
\end{equation}
In the limit $\varphi \rightarrow 0$ or
$\varphi \rightarrow \pi $ we
have
\begin{equation}
\lim\limits_{\varphi \rightarrow \frac
\pi 2\mp \frac \pi 2}\mathcal{N}_\psi
=0  \label{a3.6}
\end{equation}
\begin{equation}
\lim\limits_{\varphi \rightarrow \frac
\pi 2\mp \frac \pi 2}\mathcal{D}_\psi
=2\pi \left( \mathcal{A}_\psi ^c\cos
\frac{\varphi _0}2+x\mathcal{A}_\psi
^s\sin \frac{\varphi _0}2\right)^2
\equiv \Lambda _{\varphi _0}(x),
\label{a3.7}
\end{equation}
where
\[
\mathcal{A}_\psi ^c=\frac{A_{\alpha ,\kappa }^{(1)}
(\frac \pi 2\mp \frac \pi
2)}{A_\alpha ^{(1)}},\qquad
\mathcal{A}_\psi ^s=\frac{\Psi (q)A_{\alpha
,\kappa }^{(2)}(\frac \pi 2\mp \frac \pi 2)}
{\Psi (q)A_\alpha ^{(2)}},
\]
and
\[
x=q^\kappa e^{i(\frac \pi 2\mp \frac \pi 2)}.
\]

It is wellknown (\cite{L11} ch.VII, sect.1$)$
that in nondefinite case
spectrum of $X_{\varphi _0}$ is discreet, real
and can have accumulation
points only in infinity. The support of the
spectral measure $\sigma
_{\varphi _0}$ coincides with the spectrum of
$X_{\varphi _0},$ and, as
follows from (\ref{a3.6}) and (\ref{a3.7}), it
coincides also with set $\Pi
_{\varphi _0}$ of zeros for function
$\Lambda _{\varphi _0}(x).$ Note that
the function $\Lambda _{\varphi _0}(x)$
is even, and to find the set
\begin{equation}
\Pi _{\varphi _0}=
\left\{ \pm x_k\right\} _{k=1}^{M\leq \infty },\qquad
0<x_1<x_2<\ldots <x_M,  \label{a3.8}
\end{equation}
we must solve the equation
\begin{equation}
\mathcal{A}_\psi ^c(0)\cos \frac{\varphi _0}2+
x\mathcal{A}_\psi ^s(0)\sin
\frac{\varphi _0}2=0,  \label{a3.9}
\end{equation}
where

\begin{equation}
\mathcal{A}_\psi ^c(0)=\frac{A_{\alpha ,\kappa }^{(1)}
(0)}{A_\alpha ^{(1)}}
,\qquad \mathcal{A}_\psi ^s(0)=
\frac{\Psi (q)A_{\alpha ,\kappa }^{(2)}(0)}{
\Psi (q)A_\alpha ^{(2)}}.  \label{a3.9a}
\end{equation}

As follows from (\ref{a2.2}) and (\ref{a2.3}),
we have (with $x=q^\kappa )$
\begin{equation}
\begin{array}{ll}
A_{\alpha ,\kappa }^{(1)}(0)=
& 1+\sum\limits_{k=1}^\infty \frac 1{\left[
2k-1\right] !!}\sum\limits_{m=0}^{k-1}
(-1)^{k-m}x^{2(k-m)}\alpha
_{2m-1,2k-2}= \\
& P_0(x)+\sum\limits_{k=1}^\infty
(-1)^k\sqrt{\frac{\left[ 2k-2\right] !!}{
\left[ 2k-1\right] !!}}xP_{2k-1}(x)
\end{array}
\label{a3.10}
\end{equation}
\begin{equation}
\begin{array}{ll}
\Psi (q)A_{\alpha ,\kappa }^{(2)}(0)=
& \Psi (q)+\sum\limits_{k=1}^\infty
\frac 1{\left[ 2k\right]
!!}\sum\limits_{m=0}^{k-1}(-1)^{k-m}
x^{2(k-m)}\alpha _{2m-1,2k-1} \\
& P_0(x)+\sum\limits_{k=1}^\infty
(-1)^k\sqrt{\frac{\left[ 2k-1\right] !!}{
\left[ 2k\right] !!}}P_{2k}(x)
\end{array}
\label{a3.11}
\end{equation}
Note that with help of the recurrent relation (\ref{i6})
we can rewrite the $n$-th partial sum for the series in
(\ref{a3.10}) as
\begin{equation}
P_0(x)+\sum\limits_{k=1}^n(-1)^k\sqrt{
\ds{\frac{\left[ 2k-2\right] !!}{\left[ 2k-1\right] !!}}}
\ xP_{2k-1}(x)=(-1)^n\sqrt{
\ds{\frac{\left[ 2n\right] !!}{\left[ 2n-1\right] !!}}}
\ P_{2n}(x)  \label{a3.12}
\end{equation}

In what follows we for simplicity we not consider the
general case of
arbitrary self adjoint extension $X_{\varphi _0}$ of
coordinate operator $X$
and restrict ourself to
the concrete cases $X_0$ ($\varphi _0=0$) and $X_\pi
$ ($\varphi _0=\pi $). In the case $\varphi _0=0$ the
support of the
spectral measure, that is the set $\Pi _0=
\left\{ \pm x_k\right\}
_{k=1}^{M\leq \infty }$
of zeros must be finded from the equation
\begin{equation}
P_0(x)+\sum\limits_{k=1}^\infty (-1)^k
\sqrt{\frac{\left[ 2k-2\right] !!}{
\left[ 2k-1\right] !!}}xP_{2k-1}(x)=0  \label{a3.13}
\end{equation}
As follows from relations (\ref{i15}),
(\ref{a3.5}), (\ref{a3.3}) and (\ref{a3.1}), the function
$\psi (\gamma ;\tau )$ for the case $\varphi _0=0$ is
given by
\begin{equation}
\psi (\gamma ;\tau )=\frac{q^\kappa }\pi
\frac{\mathrm{Im}\left( \mathcal{B}
_{\alpha ,\beta ,\kappa }^{(-)}+
\mathcal{B}_{\alpha ,\beta ,\kappa
}^{(+)}\right) }{\mathcal{A}_{\alpha ,\kappa }^{(-)}+
\mathcal{A}_{\alpha
,\kappa }^{(+)}}.  \label{a3.14}
\end{equation}

Analogously for the case $\varphi _0=\pi $
the set of zeros $\Pi _\pi
=\left\{ \pm x_k\right\} _{k=1}^{M\leq \infty }$
consist from the roots of
the equation
\begin{equation}
x(P_0(x)+\sum\limits_{k=1}^\infty
(-1)^k\sqrt{\frac{\left[ 2k-1\right] !!}{
\left[ 2k\right] !!}}P_{2k}(x))=0,  \label{a3.16}
\end{equation}
and function $\psi (\gamma ;\tau )$ is given by
\begin{equation}
\psi (\gamma ;\tau )=
\frac{q^\kappa }\pi \frac{\mathrm{Im}\left( \mathcal{B}
_{\alpha ,\beta ,\kappa }^{(+)}-
\mathcal{B}_{\alpha ,\beta ,\kappa
}^{(-)}\right) }{\mathcal{A}_{\alpha ,\kappa }^{(+)}-
\mathcal{A}_{\alpha ,\kappa }^{(-)}}.  \label{a3.17}
\end{equation}
\section{Construction of the spectral measures
$\sigma _0$ and $\sigma _\pi . $}
\setcounter{equation}{0}

\subsection{The case of the measure
$\sigma _0$ ($\varphi _0=0$)}

From (\ref{a3.6}), it follows that
$\sigma _0(\Delta )=0$ if interval $
\Delta $ does not contain points $x_k\in \Pi _0$ and
in the case when $
\Delta $ contains only the single point $x_k,$ we have
\begin{equation}
\sigma _0(\Delta )\equiv \sigma _0(x_k)>0.  \label{a4.1}
\end{equation}
According to (\ref{i15}) and (\ref{a4.1}),
for $x_k\in \Pi _0$ ( we take for
definiteness $x_k>0$)
\begin{equation}
\sigma _0(x_k)=
\lim\limits_{\delta \rightarrow 0}\lim\limits_{\tau
\rightarrow 0}\int\limits_{{x_k}-\delta }^{{x_k}+
\delta }\psi
(\gamma ;\tau )d\gamma  \label{a4.2} \end{equation} For
computation of the right hand side of the relation (\ref{a4.2})
we decompose function $\psi (\gamma ;\tau )$ in the vicinity of
the point $ z=\gamma +i\tau =(x_k,0).$ Using (\ref{a3.1}),
(\ref{a3.3}) and (\ref{a3.14} ) we obtain \begin{equation} \psi
(\gamma ;\tau )=\frac{q^\kappa }\pi \frac{\mathrm{Im}\left(
e^{i\varphi }\overline{A_{\alpha ,\kappa }^{(1)}(\varphi )}\Phi
(q)A_{\beta ,\kappa }^{(2)}(\varphi )\right) }{\left( A_{\alpha
,\kappa }^{(1)}\right) ^2}
\label{a4.3}
\end{equation}
Using the relations
\[
q^\kappa e^{i\varphi }=\gamma +i\tau \Leftrightarrow
q^\kappa =\sqrt{\gamma ^2+\tau ^2},\quad e^{i\varphi
}=\frac{\gamma +i\tau }{\sqrt{\gamma ^2+\tau ^2 }}
\]
and taking into account the notation (\ref{a2.10})
we rewrite (\ref{a4.3}) in the form
\begin{equation}
\psi (\gamma ;\tau )=\frac 1\pi \frac{\gamma G_\gamma +
\tau G_\tau }{\left(
A_{\alpha ,\kappa }^{(1)R}\right) +
\left( A_{\alpha ,\kappa }^{(1)I}\right) }
,  \label{a4.4}
\end{equation}
where
\[
\begin{array}{ll}
G_\gamma = & A_{\alpha ,\kappa }^{(1)R}\Phi (q)
A_{\beta ,\kappa }^{(2)I}-A_{\alpha ,\kappa }^{(1)I}
\Phi (q)A_{\beta ,\kappa }^{(2)R} \\[5pt]
G_\tau = & A_{\alpha ,\kappa }^{(1)R}\Phi (q)
A_{\beta ,\kappa }^{(2)R}+A_{\alpha ,\kappa }^{(1)I}
\Phi (q)A_{\beta ,\kappa }^{(2)I}
\end{array}
\]
Using the relations
\begin{equation}
\begin{array}{ll}
\mathrm{Re}\left( \gamma +i\tau \right) ^{2s} &
=\gamma ^{2s}-C_{2s}^2\gamma
^{2s-2}\tau ^2+\ldots +(-1)^s\tau ^{2s} \\[5pt]
\mathrm{Im}\left( \gamma +i\tau \right) ^{2s} &
=\gamma \tau \left\{
C_{2s}^1\gamma ^{2s-2}-C_{2s}^3\gamma ^{2s-4}\tau ^2+\ldots
+(-1)^{s-1}C_{2s}^{2s-1}\tau ^{2s-2}\right\}
\end{array}
\label{a4.5}
\end{equation}
and formulas (\ref{a2.2}),
(\ref{a2.5}) and (\ref{a2.10}) we obtain
\begin{equation}
\begin{array}{ll}
A_{\alpha ,\kappa }^{(1)R}
(\gamma ,\tau ) & =1+\sum\limits_{k=1}^\infty
\frac 1{\left[ 2k-1\right] !!}
\sum\limits_{m=0}^{k-1}(-1)^{k-m}\left\{
\gamma ^{2(k-m)}-\right. \\[5pt]
& \left. -C_{2(k-m)}^2\gamma ^{2(k-m-1)}
\tau ^2+\ldots +(-1)^{k-m}\tau
^{2(k-m)}\right\} \alpha _{2m-1,2k-2};
\end{array}
\label{a4.6}
\end{equation}
\begin{equation}
\begin{array}{ll}
A_{\alpha ,\kappa }^{(1)I}(\gamma ,\tau ) & =\gamma \tau
\sum\limits_{k=1}^\infty \frac 1{\left[ 2k-1\right]
!!}\sum\limits_{m=0}^{k-1}(-1)^{k-m}\left\{ C_{2(k-m)}^1\gamma
^{2(k-m-1)}-\right. \\[5pt]
& \left. -C_{2(k-m)}^3\gamma ^{2(k-m-2)}\tau ^2+\ldots
+(-1)^{k-m-1}C_{2(k-m)}^{2(k-m)-1}
\tau ^{2(k-m-1)}\right\} \alpha
_{2m-1,2k-2}
\end{array}
\label{a4.7}
\end{equation}
\begin{equation}
\begin{array}{ll}
\Phi (q)A_{\beta ,\kappa }^{(2)R}(\gamma ,\tau ) &
=\Phi (q)+\sum\limits_{k=2}^\infty
\frac 1{\left[ 2k-1\right]
!!}\sum\limits_{m=0}^{k-2}(-1)^{k-m-1}
\left\{ \gamma ^{2(k-m-1)}-\right.
\\[5pt]
& \left. -C_{2(k-m-1)}^2\gamma ^{2(k-m-2)}
\tau ^2+\ldots +(-1)^{k-m-1}\tau
^{2(k-m-1)}\right\} \beta _{2m,2k-2}
\end{array}
\label{a4.8}
\end{equation}
\begin{equation}
\begin{array}{ll}
\Phi (q)A_{\beta ,\kappa }^{(2)I}(\gamma ,\tau ) &
\!=\!\gamma \tau
\sum\limits_{k=1}^\infty \frac 1{\left[ 2k-1\right]
!!}\sum\limits_{m=0}^{k-2}(-1)^{k-m-1}
\left\{ C_{2(k-m-1)}^1\gamma
^{2(k-m-1)}-\right. \\
& \left. \! -C_{2(k-m-1)}^3\gamma ^{2(k-m-2)}
\tau ^2+\!\ldots\!+(-1)^{k-m-2}C_{2(k-m-1)}^{2(k-m-1)-1}
\tau ^{2(k-m-2)}\right\} \beta _{2m,2k-2}
\end{array}
\label{a4.9}
\end{equation}
Let $\widetilde{\gamma }:=
\gamma -x_k$ and $f_1\left( \widetilde{\gamma }
;\tau \right) =$
\[
\begin{array}{ll}
f_1\left( \widetilde{\gamma };\tau \right)
=A_{\alpha ,\kappa }^{(1)R}(
\widetilde{\gamma }+x_k;\tau ) &
\varphi _1\left( \widetilde{\gamma };\tau
\right) =\Phi (q)A_{\beta ,\kappa }^{(2)R}
(\widetilde{\gamma }+x_k;\tau ) \\
f_2\left( \widetilde{\gamma };\tau \right)
=A_{\alpha ,\kappa }^{(1)I}(
\widetilde{\gamma }+x_k,\tau ) &
\varphi _2\left( \widetilde{\gamma };\tau
\right) =\Phi (q)A_{\beta ,\kappa }^{(2)I}
(\widetilde{\gamma }+x_k;\tau )
\end{array}
\]
If we expand this functions in Teylor series about
point $M(0,0)$ and save
leading terms we obtain
\begin{equation}
f_1\left( \widetilde{\gamma };\tau \right) =
\widetilde{\gamma }c_1(
\widetilde{\gamma })+
\tau ^2c_2(\widetilde{\gamma })+
\tau ^4c_3(\widetilde{
\gamma };\tau )+\ldots  \label{a4.10}
\end{equation}
where
\begin{equation}
\begin{array}{ll}
c_1(\widetilde{\gamma })=
\sum\limits_{k=1}^\infty
\ds{\frac{\partial ^kf_1}{\partial
\widetilde{\gamma }^k}}(0;0) \
\ds{\frac{\widetilde{\gamma }^{k-1}}{k!}} &
c_1(0)=\ds{\frac{\partial f_1}
{\partial \widetilde{\gamma }}}(0;0) \\[3pt]
c_2(\widetilde{\gamma })=\sum\limits_{k=0}^\infty
\ds{\frac{\partial ^{k+2}f_1}
{\partial \widetilde{\gamma }^k\partial \tau ^2}}(0;0) \
\ds{\frac{\widetilde{\gamma }^{k-1}}{(k+2)!}} &
c_2(0)=\ds{\frac{\partial ^2f_1}{\partial \tau ^2}}(0;0) \
\ds{\frac 1{2!}} \\[3pt]
c_3(\widetilde{\gamma };\tau )=
\sum\limits_{s=2}^\infty \tau^{2(s-2)}
\sum\limits_{k=0}^\infty
\ds{\frac{\partial ^{k+2s}f_1}
{\partial \widetilde{\gamma }^k\partial \tau ^{2s}}}(0;0) \
\ds{\frac{\widetilde{\gamma }^k}{ (k+2s)!}} &
c_3(0;0)=\ds{\frac{\partial ^4f_1}{\partial \tau ^4}}(0;0) \
\ds{\frac 1{4!}}
\end{array}
\label{a4.11}
\end{equation}
\begin{equation}
f_2\left( \widetilde{\gamma };\tau \right) =
\tau c_4(\widetilde{\gamma }
)+\tau ^3c_5(\widetilde{\gamma };\tau )+\ldots  \label{a4.14}
\end{equation}
\begin{equation}
\begin{array}{ll}
c_4(\widetilde{\gamma })=\sum\limits_{k=0}^\infty
\ds{\frac{\partial ^{k+1}f_2}
{\partial \widetilde{\gamma }^k\partial \tau }}(0;0) \
\ds{\frac{\widetilde{\gamma }^k}{(k+1)!}} &
c_4(0)=\ds{\frac{\partial f_2}
{\partial \widetilde{\gamma }}}(0;0) \\[3pt]
c_5(\widetilde{\gamma };\tau )=
\sum\limits_{s=1}^\infty \tau ^{2(s-1)}
\sum\limits_{k=0}^\infty
\ds{\frac{\partial ^{k+2s+1}f_2}
{\partial \widetilde{\gamma }^k\partial \tau ^{2s+1}}}(0;0) \
\ds{\frac{\widetilde{\gamma }^k}{(k+2s+1)!}} &
c_5(0;0)=\ds{\frac{\partial ^3f_2}{\partial \tau ^3}}(0;0)
\ds{\frac 1{3!}}
\end{array}
\label{a4.15}
\end{equation}
\begin{equation}
\varphi _1\left( \widetilde{\gamma };\tau \right) =
c_6(\widetilde{\gamma })+\tau ^2c_7(\widetilde{\gamma })+
\tau ^4c_8(\widetilde{\gamma };\tau)+\ldots  \label{a4.17}
\end{equation}
\begin{equation}
\begin{array}{ll}
c_6(\widetilde{\gamma })=\sum\limits_{k=0}^\infty
\ds{\frac{\partial ^k\varphi _1}
{\partial \widetilde{\gamma }^k}}(0;0) \
\ds{\frac{\widetilde{\gamma }^k}{(k)!}} &
c_6(0)=\varphi _1(0;0) \\[3pt]
c_7(\widetilde{\gamma })=
\sum\limits_{k=0}^\infty
\ds{\frac{\partial ^{k+2}\varphi _1}
{\partial \widetilde{\gamma }^k\partial \tau ^2}}(0;0) \
\ds{\frac{\widetilde{\gamma }^k}{(k+2)!}} &
c_7(0)=\ds{\frac{\partial ^2\varphi _1}
{\partial \tau ^2}}(0;0)
\ds{\frac 1{2!}} \\[5pt]
c_8(\widetilde{\gamma };\tau )=
\sum\limits_{s=2}^\infty \tau^{2(s-2)}
\sum\limits_{k=0}^\infty
\ds{\frac{\partial ^{k+2s}\varphi _1}
{\partial \widetilde{\gamma }^k\partial \tau ^{2s}}}(0;0) \
\ds{\frac{\widetilde{\gamma }^k}{(k+2s)!}} &
c_8(0)=\ds{\frac{\partial ^4\varphi _1}
{\partial \tau ^4}}(0;0) \ \ds{\frac 1{4!}}
\end{array}
\label{a4.18}
\end{equation}
\begin{equation}
\left( \widetilde{\gamma };\tau \right) =
\tau c_9(\widetilde{\gamma })+
\tau ^3c_{10}(\widetilde{\gamma };\tau )+
\ldots  \label{a4.21}
\end{equation}
\begin{equation}
\begin{array}{ll}
c_9(\widetilde{\gamma })=
\sum\limits_{k=0}^\infty
\ds{\frac{\partial ^{k+1}\varphi _2}
{\partial \widetilde{\gamma }^k\partial \tau }}(0;0)\
\ds{\frac{\widetilde{\gamma }^k}{(k+1)!}} & c_9(0)=
\ds{\frac{\partial \varphi _2}{\partial \tau }}(0;0)
\\[5pt]
c_{10}(\widetilde{\gamma };\tau )=
\sum\limits_{s=1}^\infty \tau ^{2(s-1)}
\sum\limits_{k=0}^\infty
\ds{\frac{\partial ^{k+2s+1}\varphi _2}
{\partial \widetilde{\gamma }^k\partial \tau ^{2s+1}}}(0;0) \
\ds{\frac{\widetilde{\gamma }^k}{(k+2s+1)!}} & c_{10}(0)=
\ds{\frac{\partial ^3\varphi _2}{\partial \tau ^3}}(0;0) \
\ds{\frac 1{3!}}
\end{array}
\label{a4.22}
\end{equation}
Using (\ref{a4.4}),
we can rewrite (\ref{a4.2}) for $x_k>0$ in the form
\begin{equation}
\begin{array}{ll}
\sigma _0(x_k)= &
\ds{\frac 1\pi} \lim\limits_{\delta \rightarrow 0}
\lim\limits_{\tau \rightarrow 0}
\int\limits_{-\delta }^\delta
\ds{\frac{ \left( \widetilde{\gamma }+x_k\right)
\left( f_1\varphi _2-\varphi _1f_2\right)
+\tau \left( f_1\varphi _1+\varphi _2f_2\right) }
{\left( f_1^{\;2}+
f_2^{\;2}\right) }}d\widetilde{\gamma } \\[5pt]
& =\ds{\frac 1\pi} \lim\limits_{\delta \rightarrow 0}
\lim\limits_{\tau \rightarrow 0}
\int\limits_{-\delta }^\delta
\ds{\frac{B_1\left( \widetilde{ \gamma };\tau \right) }
{B_2\left( \widetilde{\gamma };\tau \right) }} \
d\widetilde{\gamma }
\end{array}
\label{a4.24}
\end{equation}
Using the formulas (\ref{a4.10}) - (\ref{a4.22})
we can rewrite quantities
$B_1\left( \widetilde{\gamma };\tau \right) $ and
$B_2\left( \widetilde{\gamma };\tau \right) $
entering (\ref{a4.24}) in the form
\begin{equation}
B_1\left( \widetilde{\gamma };\tau \right) =
\tau \left( b_1(\widetilde{\gamma })+
\widetilde{\gamma }b_2(\widetilde{\gamma })\right) +
\tau ^3b_3(\widetilde{\gamma };\tau )
\label{a4.25}
\end{equation}

\begin{equation}
\begin{array}{ll}
B_2\left( \widetilde{\gamma };\tau \right) & =
\widetilde{\gamma }c_1^{\;2}(\widetilde{\gamma })+
\tau ^2b_4(\widetilde{\gamma })+
\tau ^4b_5(\widetilde{\gamma };\tau ) \\[5pt]
& =\alpha (\widetilde{\gamma };\tau )+
\tau ^4b_5(\widetilde{\gamma };\tau )
\end{array}
\label{a4.27}
\end{equation}
where
\begin{equation}
\begin{array}{l}
b_1(\widetilde{\gamma })=
-c_4(\widetilde{\gamma })c_6(\widetilde{\gamma })
\\[5pt]
b_2(\widetilde{\gamma })=
\left( \widetilde{\gamma }+x_k\right)
\left( c_1(\widetilde{\gamma })
c_9(\widetilde{\gamma })- c_4(\widetilde{\gamma })
c_6(\widetilde{\gamma })\right) +
c_1(\widetilde{\gamma })c_6(\widetilde{\gamma })
\\[5pt]
b_4(\widetilde{\gamma })=2\widetilde{\gamma }
c_1(\widetilde{\gamma })c_2(\widetilde{\gamma })+
c_4^{\;2}(\widetilde{\gamma }),\qquad b_4(0)=c_4^{\;2}(0).
\end{array}
\label{a4.26}
\end{equation}
In what follows we do not use the expressions for
$b_3(\widetilde{\gamma };\tau )$ and
$b_5(\widetilde{\gamma };\tau )$ and note only that
\begin{equation}
\begin{array}{ll}
b_3(0;\tau )= & x_k\left\{ (c_2(0)+\tau ^2c_3(0))(c_9(0)
+\tau^2c_{10}(0))-\right. \\[5pt]
& \left. (c_7(0)+\tau ^2c_8(0))(c_4(0)+
\tau ^2c_5(0))\right\} + \\[5pt]
& (c_2(0)+\tau ^2c_3(0))(c_6(0)+
\tau ^2c_7(0)+\tau ^4c_8(0))+ \\[5pt]
& (c_4(0)+\tau ^2c_5(0))(c_9(0)+\tau ^2c_{10}(0))
\end{array}
\label{a4.29}
\end{equation}
\begin{equation}
b_5(0;0)=c_2^{\;2}(0)+2c_4(0)c_5(0)  \label{a4.30}
\end{equation}
From (\ref{a4.11}), definition of $f_1$ and (\ref{a4.6})
we have
\begin{equation}
\begin{array}{ll}
c_1(0)= & \sum\limits_{l=1}^\infty
\ds{\frac 1{\left[ 2l-1\right]!!}}
\sum\limits_{m=0}^{l-1}(-1)^{l-m}2
\left( l-m\right) x_k^{2(l-m)-1}\alpha
_{2m-1,2l-2}= \\[5pt] &  \ds{\frac d{dx}}
\left( P_0(x)+\sum\limits_{l=1}^\infty (-1)^l
\sqrt{\ds{\frac{\left[ 2l-2\right] !!}
{\left[ 2l-1\right] !!}}}
xP_{2l-1}(x)\right) (x_k)
\end{array}
\label{a4.31}
\end{equation}
Analogously from (\ref{a4.18}) and (\ref{a4.8}) we have
\begin{equation}
\begin{array}{ll}
c_6(0)= & \Phi (q)+\sum\limits_{l=2}^\infty
\ds{\frac 1{\left[ 2l-1\right]!!}}
\sum\limits_{m=0}^{l-2}(-1)^{l-m-1}
x_k^{2(l-m-1)}\beta _{2m,2l-2}=
\\[5pt]
& Q_1(x_k)+\sum\limits_{l=2}^\infty (-1)^{l-1}
\sqrt{\ds{\frac{\left[ 2l-2\right]!!}
{\left[ 2l-1\right] !!}}}
Q_{2l-1}(x_k)
\end{array}
\label{a4.32}
\end{equation}
Thus from (\ref{a4.24}) and (\ref{a4.25}) we obtain
\begin{equation}
\begin{array}{ll}
\sigma _0(x_k)= &
\ds{\frac 1\pi} \lim\limits_{\delta \rightarrow 0}
\lim\limits_{\tau \rightarrow 0}\tau
\left\{ \int\limits_{-\delta }^\delta
\ds{\frac{x_kb_1\left( \widetilde{\gamma }\right) }
{B_2\left(\widetilde{\gamma };\tau \right) }} \ d \
\widetilde{\gamma }+\int\limits_{-\delta}^\delta
\ds{\frac{\widetilde{\gamma }b_2
\left( \widetilde{\gamma }\right) }
{ B_2\left( \widetilde{\gamma };\tau \right) }}\ d \
\widetilde{\gamma }+\tau ^2\int\limits_{-\delta }^\delta
\ds{\frac{b_3\left( \widetilde{\gamma }\right) }
{B_2\left( \widetilde{\gamma };\tau \right) }} \ d\
\widetilde{\gamma } \right\}
\end{array} \label{a4.33}
\end{equation}
Using (\ref{a4.27}), (\ref{a4.29}) and relations
\[
B_2\left( \widetilde{\gamma };\tau \right) =
\ds{\frac 1{\alpha \left( \widetilde{ \gamma };\tau \right) }}
\left( 1-\tau ^4
\ds{\frac{b_5\left( \widetilde{\gamma };\tau \right) }
{\alpha \left( \widetilde{\gamma };\tau \right) }}
+\ldots \right)
\]
\begin{equation}
c_4(0)=c_1(0)+2x_k\sum\limits_{l=1}^\infty
\ds{ \frac{\left[ 2l-2\right] !!}{ \left[ 2l-1\right] !!}}
\Psi _{l-1}(q)  \label{a4.34}
\end{equation}
we see that in the limit
$\delta \rightarrow 0,$\  $\tau \rightarrow 0$
\begin{equation}
\begin{array}{ll}
\tau \int\limits_{-\delta }^\delta
\ds{\frac{b_3\left( \widetilde{\gamma }\right) }
{B_2\left( \widetilde{\gamma };\tau \right) }}
\ d \widetilde{\gamma }
\sim & b_3\left( 0;\tau \right) \tau
\int\limits_{-\delta }^\delta
\ds{\frac{d \widetilde{\gamma }}
{c_1^{\;2}(0)\widetilde{\gamma }^2+
\tau ^2c_4^{\;2}(0)}} =
\ds{\frac{2b_3\left( 0;\tau \right) }{c_1(0)c_4(0)}}
\arctan \ds{\frac{\delta c_1(0)}{ \tau c_4(0)]}}.
\end{array}
\label{a4.35}
\end{equation}
So in the limit $\tau \rightarrow 0$
the 3-rd term in the sum in (\ref{a4.33} ) disappears.
The same is truth in the case when $c_1(0)=0.$

Quite similarly when
$\delta \rightarrow 0,$ $\tau \rightarrow 0$
\begin{equation}
\tau \int\limits_{-\delta }^\delta
\frac{\widetilde{\gamma }b_2\left(
\widetilde{\gamma }\right) }
{B_2\left( \widetilde{\gamma };\tau \right) }d
\widetilde{\gamma }\sim \tau b_2\left( 0\right)
\int\limits_{-\delta
}^\delta \frac{\widetilde{\gamma }d
\widetilde{\gamma }}{c_1^{\;2}(0)
\widetilde{\gamma }^2+\tau ^2c_4^{\;2}(0)}=0
\label{a4.36}
\end{equation}
Note, that (\ref{a4.36}) still hold
also in the case when $c_1(0)=0.$

{Finally, in the limit
$\delta \rightarrow 0,$ $\tau \rightarrow 0
$}

\begin{equation}
\tau \int\limits_{-\delta }^\delta
\frac{x_kb_1\left( \widetilde{\gamma }
\right) }{B_2\left( \widetilde{\gamma };
\tau \right) }d\widetilde{\gamma }
\sim -\frac{2x_kc_6\left( 0\right) }
{c_1\left( 0\right) }\arctan \frac{
\delta c_1(0)}{\tau c_4(0)}  \label{a4.37}
\end{equation}
From (\ref{a4.33}), (\ref{a4.35}) - (\ref{a4.37}) we have
\begin{equation}
\sigma _0(x_k)=\left\{
\begin{array}{lll}
0 & \mathrm{if} & c_1(0)=0, \\
-x_k \left( \frac{c_6\left( 0\right) }{c_1\left( 0\right) }
\right) & \mathrm{if} & c_1(0)\neq 0.
\end{array}
\right.  \label{a4.38}
\end{equation}
According to (\ref{a4.31}), (\ref{a4.32}) and (\ref{a4.38})
we have $\sigma
_0(x_k)=0$ if
\[
\frac d{dx}\left( A_{\alpha ,\kappa }^{(1)}(0)\right)
(x_k)=\left\{
\begin{array}{l}
0 \\
\infty
\end{array}
\right.
\]
otherwise
\begin{equation}
\sigma _0(x_k)=
-x_k  \frac{Q_1(x_k)+\sum\limits_{l=2}^
\infty (-1)^{l-1}\sqrt{\frac{\left[ 2l-
2\right] !!}{\left[ 2l-1\right] !!}}
Q_{2l-1}(x_k)}{\frac d{dx}\left( P_0(x)+
\sum\limits_{l=1}^\infty (-1)^l\sqrt{
\frac{\left[ 2l-2\right] !!}{\left[ 2l-1\right] !!}}
xP_{2l-1}(x)\right) (x_k)
} \label{a4.39}
\end{equation}
One can to check the positiveness of the measure
$\sigma_0$ \cite{L10}.

\subsection{The case of the measure
$\sigma _\pi $ ($\varphi _0=\pi $).}

Quite similarily one can consider the case when
$\varphi _0=\pi.$ We omit
related calculations and give only final expression
for the measure. Namely, $\sigma _\pi (x_k)=0$ if
\[
\frac d{dx}\left( \Psi (q)A_{\alpha ,\kappa }^{(2)}
(0)\right) (x_k)=\left\{
\begin{array}{l}
0 \\
\infty
\end{array}
\right.
\]
otherwise
\begin{equation}
\sigma _\pi (x_k)=\frac 1{x_k}\frac{1+\sum\limits_{l=1}^\infty
\frac 1{\left[ 2l\right]!!}\sum\limits_{m=0}^{l-1}
(-1)^{l-m}x_k^{2(l-m)}\beta _{2m,2l-1}}
{\sum\limits_{l=1}^\infty \frac 1{\left[ 2l\right]!!}
\sum\limits_{m=0}^{l-1}(-1)^{l-m}2\left( l-m\right)
x_k^{2(l-m)-1}\alpha _{2m-1,2l-1}}
\label{a4.40}
\end{equation}
In view of (\ref{a1.1})-(\ref{a1.4})one can to rewrite
(\ref{a4.40}) as
\begin{equation}
\sigma _\pi (x_k)=\frac 1{x_k}\frac{
Q_1(x_k)+\sum\limits_{l=1}^\infty (-1)^l
\sqrt{\frac{\left[ 2l-1\right] !!}{
\left[ 2l\right] !!}}x_kQ_{2l}(x_k)}{\frac d{dx}\left(
P_0(x)+\sum\limits_{l=1}^\infty (-1)^l
\sqrt{\frac{\left[ 2l-1\right] !!}{
\left[ 2l\right] !!}}P_{2l}(x)\right) (x_k)} \label{a4.41}
\end{equation}

\section{Conclusion}

As follows from (\ref{i15}), (\ref{i16}), (\ref{a3.13}),
(\ref{a3.16}) and general results  from \cite{L10}
(ch.3 \S4) measure $\sigma_0$ which correspond to Stieltjes
transformation $m_0(z)$ from (\ref{i15}) is equal to
$\sigma_{\pi}$ which correspond to the selfadjoint extension
$\overline{X}$ with $\phi_0=\pi$. At the same time measure
$\sigma_{\infty}$ which correspond to $m_{\infty}(z)$ connected
 with the selfadjoint extension $\overline{X}$ with $\phi_0=0$.

It is not too hard to find general connection between values of
the parameter $t=t_0$ and the related selfadjoint extensions
$\overline{X}$ labeled by $\phi_0$ which looks as
\be\l{con1}
t_0=-{\rm ctg}\ds{\frac{\phi_0}{2}}
\ds{ \frac {\psi A_{\alpha}^{(2)}} {A_{\alpha}^{(1)}} }\, .
\ee

It will be interesting to obtain the interpretation of the
entire function
\be\l{con2}
F_t(z;q) =\left[A^{(1)}_{\alpha , \kappa}(0)\right] (z) t -
\left[\psi(q)A^{(2)}_{\alpha , \kappa}(0)\right] (z)\, ,
\ee
which roots are supports of the spectral measure $\sigma_t$,
in terms of a one or
another standard $q$-special function.

We would like to stress once more that
 all results reported in this work
 are hold for each Jacobi matrix (\ref{d6})
which entries fulfill the restrictions
\[ \begin{array}{c} \left[b_ {n}\right]
\left[b_{n+2}\right] \leq \left[b_{n+1}\right] ^2,\quad n\geq 0
\\ \sum_{n=0}^\infty \left[b_{n}\right] ^{-\frac 12}<\infty
\end{array}
\]
which provide the associated moment problem to be indetermine.

This allows to expect that our results may be useful also
in other problems in which indetermine moment problem arised.


\newpage

\begin{thebibliography}{999}\small

\bibitem{L1} Coon D.D., Baker M.,
        {\it Phys. Rev. D.}, {\bf D2}, 2349 (1970); \\
        Coon D.D., Yu S., Baker M.M.,
        {\it Phys. Rev. D.}, {\bf D5}, 1429 (1972);

\bibitem{L2}
{\rm Arik M., Coon D.D.},
{\it Hilbert space of analytical functions and
generalised coherent states.}
{\it J.  Math. Phys.}, {\bf 17}, no.4, 524-527 (1976);

\bibitem{L3}
{\rm Kuryskin V.V.},
Manuscript deponented in VINITI, no 3937-76, 1976 (in Russian)
\\
{\rm Kuryskin V.V.},
{\it Ann. Found. L. de Broigle}, {\bf 5},no~2, 111-126 (1980)

\bibitem{L4}
Biedenharne L.C.,
{\it The quantum group $su_q(2)$ and
$q$-analog of the boson operators},
{\it J. Phys. A.}, {\bf 22}, no~18 (1989) L 873-878

\bibitem{L5}
Macfarlane A.J.,
{\it On $q$-analogues of the quantum harmonic
oscillator and the quantum group $su_q(2)$}
{\it J. Phys. A.}, {\bf 22}, no~21 (1989) 4581-4586

\bibitem{L6}
Kulish P.P., Damaskinsky E.V.,
      {\it J. Phys. A.}, {\bf 23}, no 9, L~415-419 (1990)

\bibitem{L7} Floreannini R.,
 {\it $q$-orthogonal polynomials and the
 $q$-oscillator quantum group,}
 {\it Lett. Math. Phys.}, {\bf 22}, no~1 (1991) 45-54;
\\
 Floreanini R., Le~Teurneux J., and Vinet L.,
 {\it More on the $q$-oscillator algebra and
  $q$-orthogonal  polynomials, }
 {\it J. Phys. A.}, {\bf 28}, no~10 (1994) L287-294;
\\
 Floreanini R.,
 {\it Quantum algebras and $q$-special functions,}
 {\it Ann. Phys.} {\bf 221} (1993).

\bibitem{L8}
 Atakishiev N.M., and Suslov S.K.,
 {\it On a realization of the $q$-harmonic oscillator,}
 {\it Teor. Mat. Fiz.},  (1990) 64;
 [{\it Theor. Math. Phys.} {\bf 85} (1990) 1055];
\\
 Atakishiev N.M., and Feinsilver P.,
 {\it On the coherent states for the $q$-Hermitean polynomial
 and related Fourier transformation},
 {\it J. Phys. A.}, {\bf 29}, no~8 (1996) 1659-1664;
\\
 Atakishiev N.M., Frank A., and Wolf K.B.,
 {\it A simple difference realization of Heisenberg $q$-algebra}
 {\it J.Math.Phys.}, {\bf 35}, no~7 (1994) 3253-3260.

\bibitem{Z199}
 Damaskinsky E.V., and P.P.Kulish P.P.,
 {\it $q$-Hermite polynomials and $q$-os\-cil\-la\-tors,\/}
 {\it Zap. Nauch. Sem. POMI,\/} {\bf 199} (1992) 81-90;
 (in Russian);

\bibitem{Z189}
Damaskinsky E.V., Kulish P.P.,
{\it Deformed  oscillators and their applications},
{\it Zap. Nauch. Seminarov LOMI}, {\bf 189}, 37-74 (1991) (in
Russian) \\ English transl:
{\it J.Soviet.Math.} {\bf 62}, 2963 (1992)

\bibitem{PPK}
Kulish P.P.,
{\it Irreducible representations of deformed
oscillator and coherent states},
Preprint KTH-96/21, 11pp.,  Stockholm, 1996.

\bibitem{L9}
Damaskinsky E.V., and Kulish P.P.,
{\it Irreducible representations of
deformed oscillator algebra and $q$-special functions,}
{\it Intern. J. Mod. Phys.A.},
1990,12,N~1,153-158; q-alg 9610002,

\bibitem{BDY}
{Borzov V.V., Damaskinsky E.V., and Yegorov S.B.},
{\it Some remarks on the representations of
the generalized deformed oscillator algebra,\/}
{\it Zap.Nauch.Seminarov LOMI} {\bf 245}, 80-106 (1997) (in
Russian) Stockholm pr-t TRITA-MAT-1995-MA-20 ({q-alg/9509022})


\bibitem{L10}
 Akhiezer N.I., {\it The classical moment problem and some
related questions in analysis,}
Hafner Publ. Co, New York 1965;

\bibitem{L11}
Berezanski{\u \i} Yu.M.,
{\it Expansions in eigenfunctions of selfadjoint operators},
{\it Transl. Math. Monographs 17}, Amer. Math. Soc.,
Providence, R.I., 1968

\bibitem{L12}
Shohat~J.A. and Tamarkin~J.D.,
{\it The problem of moments},
Amer. Math. Soc., Providence, R.I., 1963

\bibitem{Chihara}
Chihara T.S.,
{\it Indetermine Hamburger moment problem},
{\it Pacif. J. Math.}, {\bf 27}, no~3 (1968) 475-484;
\\
Chihara T.S.,
{\it Indeterminate symmetric moment problems},
{\it J. Math. Anal. Appl.}, {\bf 85}, 331-346 (1982).

\bibitem{Chihara2}
Chihara T.S.,
{\it An Introduction to Orthogonal Polynomials}
 Math. and Appl. 13,  Gordon and Breach, New York, 1978

\bibitem{Nagel}
Nagel B.,
{\it Higher power squeezed states, Jacobi matrices, and
the Hamburger moment problem},
Contribution to the 1997
Balatonfured conference on squeezed states.

\bibitem{L13}
Burban I.M., and Klymik A.U.,
{\it On spectral properties of $q$-oscillator operators}
{\it Lett. Math. Phys.}, {\bf 29} (1993) 13-18
[Pr-t ITP-92-59E, Kiev].

\bibitem{L14}
Chung W.-S., and  Klimyk A.U.,
{\it On position and momentum operators in the $q$-oscillator
algebra},
{\it J. Math. Phys.}, {\bf 37}, no~2 (1996) 917-932

\bibitem{Ex}
Exton H., {\it q-hypergeometric functions and  applications},
Chichester. Ellis Horwood.  1983;

\bibitem{Gasper}
{Gasper G., and Rahman M.},
{\it Basic Hypergeometric Series}.
{Cambridge Univ.Press}, 1990 in
{\it Encyclopedia of Mathematics and its
Applications 35}, Cambridge, 1990.

\bibitem{Andrews}
Andrews G.A.,
{\it $q$-Series, their development and
application in analysis, number theory, combinatorics,
physics and computer mathematics}.
AMS regional conference series 66, 1986.


\bibitem{b1} Moak D.S.
{\it The $q$-analogue of the Laguerre polynomials},
{\it J. Math. Anal. Appl.}, {\bf 81} 20-47 (1981)

\bibitem{b4}
Berg C.,
{\it Indetermine moment problem},
{\it J. Comp. Appl. Math.} {\bf 65}, 27-55 (1965)
\\
Berg C.,  Duran A.J.,
{\it Orthogonal polynomials, $L^2$ spaces and
entire functions},
{\it Math.Scand.} {\bf 79}, 209-223 (1996)

\bibitem{L102}
Akhiezer N.I.,
{\it Uspechi Matem. Nauk.} no.9 (1941)

\bibitem{a12}
Plesner A.I.,
{\it Uspechi Matem. Nauk.} no.9 (1941)

\bibitem{CK}
Chaichian M., and Kulish P.P.,
{\it Phys.Lett.B.,} {\bf 234},  no.1/2, 72-80 (1990);
Pr-t. CERN-TH 5969/90,

\bibitem{PK}
Kulish P.P.,
{\it Teor. Math. Phys.} {\bf 85}, no 1, 158-161(1991)

\bibitem{FdP}
Rampacher H., Stumpf H., and Wagner F.,
{\it Fortschr. Phys.}, {\bf 13}, 385 (1965)

\bibitem{L16}  Curtright T.L.,
 in: {\it Quantum Groups},
(T.L.~Curtright, D.~Fairlie and C.K.~Zachos (Eds.)),
Proceedings of the Argonne Workshop
(World Scientific, Singapore, 1990).

\bibitem{KQ} Katriel J., and Quesne C.,
{\it J. Math. Phys.}, {\bf 37}, no~4, 1650-1661 (1996)

\bibitem{h1} Rogers L.J.,
{\it On the expansion of some infinite products},
{\it Proc. London Math. Soc.}, {\bf 24}, 337-352 (1893);
{\it Second memoir on the expansion of some infinite products},
{\it Proc. London Math. Soc.}, {\bf 25}, 319-343 (1894);
{\it Third memoir on the expansion of some infinite products},
{\it Proc. London Math. Soc.}, {\bf 26}, 15-32   (1895);

\bibitem{h2} Szeg{\" o} G.,
{\it Ein Beitrag zur Theorie der Thetafunctionen},
{\it Sitz. Preuss. Ak. Wiss.} Phys.-Math.Class.,1926, 242-251;
\\
Carlitz L.,
{\it Some polynomials related to theta functions},
{\it Annali di Matematica Pura et Applicata}, {\bf 41}, 359-371
(1955);
\\
Carlitz L.,
{\it Some polynomials related to theta functions},
{\it Duke Math. J.}, {\bf 24} 521-527 (1957)

\bibitem{h3} Al-Salam W.A., and Carlitz, L.,
{\it Some orthogonal $q$-polynomials},
{\it Math. Nachr,}, {\bf 30}, no~1, 47-61 (1965)

\bibitem{h4} Cigler J.,
{\it Operatormethoden f{\" u}r $q$-Identit{\" a}en},
{\it Mn. Math.}, {\bf 88}, no~1, 87-105 (1979);

\bibitem{h5} Allaway W.R.,
{\it Some properties of the $q$-Hermite polynomials},
{\it Canadian J. Math.}, {\bf 32}, no~3, 686-694 (1980);


\bibitem{h7}  Mir - Kasimov R.M., Preprint IC/90/383, 1990

\end{thebibliography}
\end{document}